\documentclass[12pt]{article}
\usepackage{amsmath,amssymb}
\usepackage{graphicx}
\newcommand{\eqdef}{\stackrel{\text{def}}{=}}

\newcommand{\n}{\nonumber \\}
\newcommand{\bm}{\boldsymbol}
\newcommand{\ignore}[1]{}
\numberwithin{equation}{section}
\newcommand{\Romannumeral}[1]{\uppercase\expandafter{\romannumeral#1}}

\newtheorem{theo}{\bf Theorem}[section]

\newtheorem{rema}[theo]{\bf Remark}

\newtheorem{prop}[theo]{\bf Proposition}
\newtheorem{defi}[theo]{\bf Definition}
\newcommand{\ma}{\hspace{0pt}}
\newcommand{\cX}{\mathcal{X}}

\allowdisplaybreaks[4]

\begin{document}

\baselineskip=20pt
\newcommand{\preprint}{
\vspace*{-20mm}\begin{flushleft}\end{flushleft}
}
\newcommand{\Title}[1]{{\baselineskip=26pt
  \begin{center} \Large \bf #1 \\ \ \\ \end{center}}}
\newcommand{\Author}{\begin{center}
  \large \bf 
  Ryu Sasaki${}$ \end{center}}
\newcommand{\Address}{\begin{center}
     Department of Physics and Astronomy, Tokyo University of Science,
     Noda 278-8510, Japan
        \end{center}}
\newcommand{\Accepted}[1]{\begin{center}
  {\large \sf #1}\\ \vspace{1mm}{\small \sf Accepted for Publication}
  \end{center}}

\preprint
\thispagestyle{empty}

\Title{Multivariate Hahn polynomials \\
and difference equations}

\Author

\Address
\vspace{1cm}

\begin{abstract}
The   multivariate Hahn polynomials are constructed explicitly as the common 
eigenvectors of  a family of second order difference
operators.  They are  orthogonal with respect to the hypergeometric multinomial distribution.
The main difference operator is adopted from the  work of Karlin-McGregor in 1975.
The minor ones are the subsets of the main one containing less and less variables.
These operators commute with each other.
In contrast to the multivariate Krawtchouk and Rahman like polynomials derived recently, 
the entire multivariate Hahn polynomials are
rational functions of the system parameters.
Complete sets of multivariate Krawtchouk and Meixner polynomials are derived by limiting procedures.
\end{abstract}

%
%
\section{Introduction}
\label{sec:intro}
Here I report the construction of multivariate Hahn polynomials as the eigenpoynomials of a family
of difference operators. Various multivariate discrete orthogonal polynomials have been discussed and built
\cite{diaconis13,dunkl0,dunkl},\cite{genest0}--\cite{KarMcG},\cite{Milch},\cite{mKrawt}--\cite{zheda},
including quite a few  multivariate Hahn polynomials 
\cite{dunkl0,genest,ilxu,ilxu2,ilxu3,KarMcG,rosengren,scarabotti,tra3,xu3}.
Construction of multivariate orthogonal polynomials as the eigenpolynomials of a certain difference operator,
natural generalisation of the Askey scheme \cite{askey}, looks quite appealing.
If the difference operator is similarity equivalent to a self-adjoint operator, the orthogonality of the polynomials
is guaranteed as the eigenvalues are non-degenerate for generic parameters.
However, determination of degree one eigenpolynomials for an $n$-variable theory usually  involves a degree $n$
characteristic equation and the eigenvalues and eigenvectors are in general irrational.
In that case, determination of higher eigenvalues and eigenvectors seems almost intractable.
For the multivariate Krawtchouk, Meixner \cite{I23,mKrawt} 
and Rahman like polynomials \cite{mRah}, this difficulty is circumvented as the forms of all higher degree
eigenpolynomials are fixed by the degree one polynomials \cite{mizu}.
The situation of the multivariate Hahn polynomials seems quite opposite.  Thanks to the high symmetry
of the adopted difference operator due to Karlin-McGregor \cite{KarMcG}, 
Definition \ref{def:Ht}, the eigenvalues of all total degree $M$ 
eigenpolynomials are the same and the eigenpolynomials are closely related to those  of single variable 
Hahn polynomials. 
Combining these `known' polynomials to a complete set of orthogonal polynomials is the main task.

My guiding principle is the orthogonality measure and difference equations, which has led to the
construction of the multivariate Krawtchouk, Meixner \cite{I23,mKrawt} 
and Rahman like polynomials \cite{mRah}.
In the Krawtchouk case, the binomial distribution is the orthogonality measure which goes naturally to the
multinomial distribution, as the binomial theorem \\
$(a+b)^N=\sum_{n=0}^N
\binom{N}{n}a^nb^{N-n}$ goes
to the multinomial theorem 
$\bigl(\sum_{i=1}^n a_i+b\bigr)^N
=\sum_{\bm{x}}\binom{N}{\bm{x}}\prod_{i=1}^na_i^{x_i}\cdot b^{N-|x|}$, $|x|\eqdef\sum_{i=1}^n x_i$.
The multivariate difference equation is the simple generalisation of the single variable one \cite{mKrawt}(3.11).
The strategy is essentially the same for the Hahn case.
The orthogonality measure of the single variable Hahn is the hypergeometric binomial distribution which derives from 
the hypergeometric binomial theorem $(a+b)_N=\sum_{n=0}^N\binom{N}{n}(a)_n(b)_{N-n}$,
in which $(a)_0=1$, $(a)_n=a(a+1)\cdots(a+n-1)$, $n\in\mathbb{N}_0$.
The multivariate case  is also governed by the hypergeometric multinomial theorem 
\begin{equation*}
\left(\sum_{i=1}^n a_i+b\right)_N=\sum_{\bm{x}}\binom{N}{\bm{x}}\prod_{i=1}^n(a_i)_{x_i}\cdot(b)_{N-|x|}.
\end{equation*}
As  for the difference operator, one could have chosen the naive multivariate 
generalisation of the single variable one $\widetilde{\mathcal H}_0$ \eqref{Ht0def}.
I choose the one $\widetilde{\mathcal H}_T$ \eqref{Htsum}
 in the seminal paper by Karlin and McGregor \cite{KarMcG}(3.13) which has many additional
operators on top of the naive one $\widetilde{\mathcal H}_0$. 
All the additional ones commute with the naive one and they share the same eigenpolynomials.
Due to the high symmetry of the adopted main difference operator 
$\widetilde{\mathcal H}_T$ \eqref{Htsum}, Definition \ref{def:Ht}, the construction of the
eigenpolynomials is straightforward.

Several papers are dedicated to the same subject, multivariate Hahn, 
\cite{dunkl0, genest, ilxu,ilxu2,KarMcG,rosengren,scarabotti,tra3,xu2,xu3}.
Many notions and functions look similar or related, but I have yet to decipher the exact relationship
with the earlier works.
For the bivariate case, the eigenpolynomials \eqref{sol1} are the same as those derived by Genest and Vinet
\cite{genest}
with some changes of notation.
But my method of derivation is totally different from others.

This paper is organised as follows.
After introducing the difference operators, $\widetilde{\mathcal H}_T$ \eqref{Htsum},
$\widetilde{\mathcal H}_0$ \eqref{Ht0def}, $\widetilde{\mathcal H}_1$ \eqref{Ht1def} in section \ref{sec:difop}, 
the corresponding real symmetric matrices $\mathcal{H}_T$ \eqref{HTdef}, etc, are derived   in section  \ref{sec:Hmat}
through the orthogonality measure defined in section \ref{sec:orthomeas}.
Two types of eigenpolynomials are introduced in section \ref{sec:eipoly}.
After a short review of various properties of the single variable Hahn polynomials
$H_m(x)$
in section \ref{sec:SHahn}, it is shown in Proposition \ref{prop:single}
that the single variable Hahn polynomials depending on various subsets of the
variable $\bm{x}$, {\em e.g,\/} $H_m(x_1+x_3+x_8)$, are the eigenpolynomials
of the main difference operator $\widetilde{\mathcal H}_T$ \eqref{Htsum} with the same eigenvalue.
Type two eigenpolynomials depending on two consecutive variables $x_i$ and $x_{>i}$ \eqref{x>def}
are introduced in sections \ref{sec:typetwo}, \ref{sec:hieig}.
Their explicit forms are shown in Theorem \ref{theo:Pishape}.
Several important properties of the type two eigenpolynomials are displayed in Proposition \ref{prop:PiForbac}, 
Theorem \ref{theo:rec} and eigenvalue Theorem \ref{theo:eigens}.
The main results, the complete set of orthogonal eigenpolynomials, are presented in Theorem \ref{theo:main},
followed by supporting Theorems \ref{theo:ii-1}, \ref{theo:genfb} and \ref{theo:eigensR}.
Multivariate Krawtchouk and Meixner polynomials are obtained  in section \ref{sec:newmKra} and \ref{sec:newmMei} 
by certain limits from the multivariate Hahn polynomials. 
Section \ref{sec:comments} is for comments. 
Appendix A  provides the derivation of the type two eigenpolynomials. 
Appendix B  gives the proof of the eigenvalue Theorem \ref{theo:eigens}.
%
%
\section{Problem setting}
\label{sec:Problem}
The purpose of the paper is to construct multivariate Hahn polynomials as eigenpolynomials
of certain difference operators. The polynomials are defined on a compact $n$-dimensional
integer lattice called $\cX$, restricted by a positive integer $N>n\ge2$,
\begin{equation}
 \bm{x}=(x_1,\ldots,x_n)\in\mathbb{N}_0^n,\quad  |x|\eqdef \sum_{i=1}^nx_i,\quad
 \mathcal{X}\eqdef\{\bm{x}\in\mathbb{N}_0^n\ |\, |x|\le N\}.
 \label{XKdef}
\end{equation}
The space of polynomials with real coefficients is $V_N(\bm{x})$
\begin{equation}
\bm{m}=(m_1,\ldots,m_n)\in\mathbb{N}_0^n,
\quad V_N(\bm{x})\eqdef{\rm Span}\left\{\bm{x}^{\bm{m}}\mid 0\le |m| \le N\right\},
\quad \bm{x}^{\bm{m}}\eqdef\prod_{i=1}^nx_i^{m_i}.
\label{VNdef}
\end{equation}
%
%
\subsection{ Difference operators}
\label{sec:difop}
The main difference operator $\widetilde{\mathcal H}_T$ depending on 
$n+1$ positive parameters $a_j>0$, $j=1,\ldots,n$ 
and $b>0$ and acting on a polynomial $p(\bm{x})\in V_N(\bm{x})$ is defined by 
the following
\begin{defi}
\label{def:Ht}
\begin{align}
\widetilde{\mathcal H}_Tp(\bm{x})&=\sum_{j=1}^n(N-|x|)(x+a_j)\bigl(p(\bm{x})-p(\bm{x}+\bm{e}_j)\bigr)\n
                                                          &+    \sum_{j=1}^nx_j(N-|x|+b)\bigl(p(\bm{x})-p(\bm{x}-\bm{e}_j)\bigr)\n
                                       & +\sum_{j\neq k\ge1}^nx_j(x_k+a_k)\bigl(p(\bm{x})-p(\bm{x}-\bm{e}_j+\bm{e}_k)\bigr),
 \end{align}  
 in which $\bm{e}_j$ is the  $j$-th unit vector, $j=1,\ldots,n$.                                       
  \end{defi}
Let us introduce 
operators $e^{\pm\partial_j}$ acting on a polynomial $p(\bm{x})\in V_N(\bm{x})$ 
\begin{equation*}
e^{\pm\partial_j}f(\bm{x})=f(\bm{x}\pm\bm{e}_j)\,e^{\pm\partial_j},
\quad \partial_j\eqdef\frac{\partial}{\partial x_j},
\quad j=1,\ldots,n.
\end{equation*}
The decomposition of the operator  $\widetilde{\mathcal H}_T$ into two parts   $\widetilde{\mathcal H}_0$
and $\widetilde{\mathcal H}_1$ is defined by the following
\begin{defi}
\label{def:Ht0}
\begin{align}
\widetilde{\mathcal H}_T&=\widetilde{\mathcal H}_0+\widetilde{\mathcal H}_1,
\label{Htsum}\\
\widetilde{\mathcal H}_0&\eqdef \sum_{j=1}^n\bigl(B_j(\bm{x})(1-e^{\partial_j})
+D_j(\bm{x})(1-e^{-\partial_j})\bigr),
\label{Ht0def}\\
&=\sum_{j=1}^n\bigl(B_j(\bm{x})-D_j(\bm{x})e^{-\partial_j}\bigr)\bigl(1-e^{\partial_j}\bigr).
\label{Ht0fac}\\
& \qquad B_j(\bm{x})\eqdef (N-|x|)(x_j+a_j),\quad D_j(\bm{x})\eqdef x_j(N-|x|+b),
\label{BDdef}\\
\widetilde{\mathcal H}_1&\eqdef \sum_{j\neq k\ge1}^n x_j(x_k+a_k)\bigl(1-e^{-\partial_j+\partial_k}\bigr).
\label{Ht1def}
\end{align}
\end{defi}
Here $\widetilde{\mathcal H}_0$ is the naive $n$-variable generalisation of the single variable difference operator
of the Hahn polynomial \cite{KLS,os12}. The functions $B_j(\bm{x})$ and $D_j(\bm{x})$ are the birth and death rates
of the $j$-th population group in the corresponding birth and death process \cite{bdsol,mKrawt}.
The additional one $\widetilde{\mathcal H}_1$ is taken from the work of Karlin-McGregor \cite{KarMcG}(3.13).
\begin{rema}
When written by $x_0\eqdef N-|x|$, $a_0=b$, the entire difference operator look very symmetric
\begin{align*}
\widetilde{\mathcal H}_T&=\sum_{j=1}^n\Bigl(x_0(x_j+a_j)(1-e^{\partial_j})
                                              +x_j(x_0+a_0)(1-e^{-\partial_j})\Bigr)\\
                                      &\quad + \sum_{j\neq k\ge1}^n x_j(x_k+a_k)\bigl(1-e^{-\partial_j+\partial_k}\bigr),
\end{align*}
except that $e^{\pm\partial_0}$ do not appear as $x_0$ not being an independent variable.  In \cite{KarMcG}
$x_0$ is treated equally as the other variables. Karlin-McGregor discuss the stochastic process, 
generalised birth and death process, rather  than the
eigenvalue problem of difference equations.
\end{rema}
A family of operators $\widetilde{\mathcal H}_i$, $i=2,\ldots,n-1$ is defined by the following
\begin{defi}
\label{def:Hti}
\begin{align}
\widetilde{\mathcal H}_i&\eqdef \sum_{j\neq k\ge i}^n x_j(x_k+a_k)\bigl(1-e^{-\partial_j+\partial_k}\bigr),
\quad i=2,\ldots,n-1.
\label{Htidef}
\end{align}
\end{defi}
It is easy to verify the following
\begin{prop}
\label{prop:Htcomm}
The difference operators commute with each other
\begin{align*}
[\widetilde{\mathcal H}_T,\widetilde{\mathcal H}_j]=0,\quad 
[\widetilde{\mathcal H}_j,\widetilde{\mathcal H}_k]=0,\quad j,k=0,1,\ldots,n-1.
\end{align*}
\end{prop}
\begin{prop}
\label{prop:zero}
A common zero-mode of all difference operators is a constant,
\begin{equation*}
\widetilde{\mathcal H}_T1=0,\quad \widetilde{\mathcal H}_i1=0,\quad i=0,1,\ldots,n-1.
\end{equation*}
$\widetilde{\mathcal H}_i$ $i\ge2$ has many zero zero modes,
\begin{equation*}
\widetilde{\mathcal H}_ip(x_1,\ldots,x_{i-1})=0,\quad i=2,\ldots,n-1,
\quad \forall p(x_1,\ldots,x_{i-1})\in V_N(\bm{x}).
\end{equation*}
\end{prop}
A subspace of $V_N(\bm{x})$ consisting of polynomials of maximal degree $M$ is defined by
\begin{equation*}
V_M(\bm{x})\eqdef{\rm Span}\left\{\bm{x}^{\bm{m}}\mid 0\le |m| \le M\right\}\subseteq V_N(\bm{x}).
\end{equation*}
\begin{prop}
\label{Hthinv}
It is obvious that the subspace $V_M(\bm{x})$ 
is invariant under the action of $\widetilde{\mathcal H}_T$ and $\widetilde{\mathcal H}_i$, $i=0,1,\ldots,n-1$,
\begin{equation*}
\widetilde{\mathcal H}_TV_M(\bm{x})\subseteq V_M(\bm{x}),\quad
\widetilde{\mathcal H}_iV_M(\bm{x})\subseteq V_M(\bm{x}),\quad i=0,1,\ldots,n-1,
\end{equation*}
and  $\widetilde{\mathcal H}_T$ and $\widetilde{\mathcal H}_i$ have eigenpolynomials in  $V_M({\bm x})$.
\begin{align}
\widetilde{\mathcal H}_Tp_{\bm m}(\bm{x})=\mathcal{E}_T(\bm{m})p_{\bm m}(\bm{x}),\quad
\widetilde{\mathcal H}_ip_{\bm m}(\bm{x})=\mathcal{E}_i(\bm{m})p_{\bm m}(\bm{x}),\quad \bm{m}\in\cX.
\label{eigpoly}
\end{align}
\end{prop}

%
%
\subsection{Orthogonality measure and real symmetric matrices corresponding to difference operators}
\label{sec:orthomeas}
As explained in \cite{bdsol,mKrawt}, the birth and death rates determine the corresponding orthogonality measure
uniquely (up to normalisation) if the compatibility conditions \cite{mKrawt}(2.17) \cite{ilxu} are satisfied.
\begin{prop}
\label{proop:Wort}
The hypergeometric multinomial distribution depending on $\bm{a}$ and $b$ is the orthogonality measure of the
difference operator system $\widetilde{\mathcal H}_T$ and $\widetilde{\mathcal H}_i$, $i=0,\ldots,n-1$,
\begin{align}
&W(\bm{x};\bm{a},b,N)=\frac{N!}{x_1!\cdots x_n!x_0!}
\frac{\prod_{i=1}^n(a_i)_{x_i}(b)_{N-|x|}}{(|a|+b)_N},
\quad \bm{a}=(a_1,\ldots,a_n)\in\mathbb{R}_{>0}^n.
\label{WHh}
\end{align} 
\end{prop}
The compatibility  conditions \cite{mKrawt}(2.17) are trivially satisfied.
\begin{align*}
\frac{B_j(\bm{x})}{D_j(\bm{x}+\bm{e}_j)}\frac{B_k(\bm{x}+\bm{e}_j)}{D_k(\bm{x}+\bm{e}_j+\bm{e}_k)}
&=\frac{\bigl(N-|x|\bigr)(x_j+a_j)}{\bigl(b+N-|x|-1\bigr)(x_j+1)}
\frac{\bigl(N-|x|-1\bigr)(x_k+a_k)}{\bigl(b+N-|x|-2\bigr)(x_k+1)}\\
&=\frac{B_k(\bm{x})}{D_k(\bm{x}+\bm{e}_k)}
\frac{B_j(\bm{x}+\bm{e}_k)}{D_j(\bm{x}+\bm{e}_k+\bm{e}_j)},
\quad j,k=1,\ldots,n.
\end{align*}
Let us introduce $\phi_0(\bm{x})$ as the square root of the orthogonality measure
\begin{equation}
\phi_0(\bm{x})\eqdef\sqrt{W(\bm{x};\bm{a},b,N)},\quad \bm{x}\in\cX,
\label{phi0def}
\end{equation}
in which the parameter dependence is suppressed for simplicity of presentation. 
It is straightforward to verify the following
\begin{prop}
\label{BDratio}
\begin{equation*}
\frac{\phi_0(\bm{x}+\bm{e}_j)}{\phi_0(\bm{x})}=\sqrt{\frac{B_j(\bm{x})}{D_j(\bm{x}+\bm{e}_j)}},
\quad j=1,\ldots,n.
\end{equation*}
\end{prop}
See \cite{mKrawt}(2.16).

%
\subsection{From $\widetilde{\mathcal H}_T$, $\widetilde{\mathcal H}_i$ operator 
to real symmetric matrices $\mathcal{H}_T$, $\mathcal{H}_i$}
\label{sec:Hmat}
Let us introduce $\mathcal{H}_T$ and $\mathcal{H}_i$ by a similarity transformation of 
$\widetilde{\mathcal H}_T$ and $\widetilde{\mathcal H}_i$,
\begin{defi}
\label{def:H}
\begin{align}
&\mathcal{H}_T\eqdef\phi_0(\bm{x})\widetilde{\mathcal H}_T\phi_0(\bm{x})^{-1},\qquad
\mathcal{H}_i\eqdef\phi_0(\bm{x})\widetilde{\mathcal H}_i\phi_0(\bm{x})^{-1},
\n
&\Longrightarrow \mathcal{H}_T\phi_0(\bm{x})=\phi_0(\bm{x})\widetilde{\mathcal H}_T,\quad
\mathcal{H}_i\phi_0(\bm{x})=\phi_0(\bm{x})\widetilde{\mathcal H}_i,\qquad i=0,1,\ldots,n-1.
\label{intert}
\end{align}
\end{defi}
It is straightforward to derive the following explicit forms of $\mathcal{H}_T$ and $\mathcal{H}_i$,
\begin{prop}
\label{prop:H}
\begin{align}
\mathcal{H}_T&=\sum_{j=1}^n\left[B_j(\bm{x})+D_j(\bm{x})-\sqrt{B_j(\bm{x})D_j(\bm{x}+\bm{e}_j)}\,e^{\partial_j}
-\sqrt{B_j(\bm{x}-\bm{e}_j)D_j(\bm{x})}\,e^{-\partial_j}\right]\n
&\ +\sum_{j\neq k}^n\left[x_j(x_k+a_k)-\sqrt{x_j(x_j+a_j-1)(x_k+1)(x_k+a_k)}\,e^{-\partial_j+\partial_k}\right]
\label{HTdef}\\
&=\sum_{j=1}^n\mathcal{A}_j^T\mathcal{A}_j+\frac12\sum_{j\neq k}^n\mathcal{B}_{j,k}^T\mathcal{B}_{j,k},
\label{semipos}\\
\mathcal{H}_i&=\frac12\sum_{j\neq k\ge i}^n\mathcal{B}_{j,k}^T\mathcal{B}_{j,k},
\qquad \qquad \qquad \qquad \qquad  i=1,\ldots,n-1,
\label{semiposi}\\
&\mathcal{A}_j\eqdef\sqrt{B_j(\bm{x})}-e^{\partial_j}\sqrt{D_j(\bm{x})},\quad
\mathcal{A}_j^T= \sqrt{B_j(\bm{x})}-\sqrt{D_j(\bm{x})}\,e^{-\partial_j},
\label{Ajdef}\\
&\mathcal{B}_{j,k}\eqdef e^{\partial_j}\sqrt{x_j(x_k+a_k)}-e^{\partial_k}\sqrt{x_k(x_j+a_j)},\n
&\mathcal{B}_{j,k}^T= \sqrt{x_j(x_k+a_k)}\,e^{-\partial_j}-\sqrt{x_k(x_j+a_j)}\,e^{-\partial_k}.
\label{Bjdef}
\end{align}
\end{prop}

These results are summarised as the following
\begin{theo}
\label{theo:HT}
$\mathcal{H}_T$ and $\mathcal{H}_i$, $i=0,\ldots,n-1$ are real symmetric positive semi-definite
matrices and the eigenvalues are real and non-negative.
In other words, $\mathcal{H}_T$ and $\mathcal{H}_i$, $i=0,\ldots,n-1$ are self-adjoint linear operators
in a $|\cX|$ dimensional compact Hilbert space $\mathbb{H}$ consisting of  real functions on $\cX$,
\begin{equation}
\mathbb{H}\eqdef \{f\Bigm| f(\bm{x})\in\mathbb{R},\ \bm{x}\in\cX\},
\label{Hildef}
\end{equation}
with the standard inner product $\langle f,g\rangle=\sum_{\bm{x}\in\cX}f(\bm{x})g(\bm{x})$.
\begin{align}
\langle \mathcal{H}_T\,f,g\rangle=\langle f,\mathcal{H}_T\,g\rangle,\quad
\langle \mathcal{H}_i\,f,g\rangle=\langle f,\mathcal{H}_i\,g\rangle,\quad i=0,1,\ldots,n-1.
\label{Hherm}
\end{align}
This also means the equalities
\begin{align}
(\widetilde{\mathcal H}_Tp,q)=(p,\widetilde{\mathcal H}_Tq),\quad
(\widetilde{\mathcal H}_ip,q)=(p,\widetilde{\mathcal H}_iq) ,\quad i=0,1,\ldots,n-1,
\label{Htherm}
\end{align}
with respect to the following inner product $(\cdot,\cdot)$ among polynomials $p,q\in V_N(\bm{x})$,
\begin{equation}
(p,q)\eqdef\langle\phi_0 p,\phi_0 q\rangle=\sum_{\bm{x}\in\cX}p(\bm{x})q(\bm{x})W(\bm{x}),\qquad
p,q\in V_N(\bm{x}).
\label{newin}
\end{equation}
The only zero-mode of $\mathcal{H}_T$ is $\phi_0(\bm{x})=\sqrt{W(\bm{x})}$.
$\mathcal{H}_T$ has a complete set of eigenvectors $\{p_{\bm m}(\bm{x})\phi_0(\bm{x})\}$,
$\bm{m}\in\cX$,  consisting of the 
eigenpolynomials of the difference operator $\widetilde{\mathcal H}_T$, \eqref{eigpoly}.
\end{theo}
$\phi_0(\bm{x})$ is the zero-mode of $\mathcal{A}_j$, $\mathcal{A}_j\phi_0(\bm{x})=0$, as is clear from 
{\bf Proposition \ref{BDratio}}.  $\phi_0(\bm{x})$ is also annihilated by $\mathcal{B}_{j,k}$,
\begin{align*}
\mathcal{B}_{j,k}\phi_0(\bm{x})&
=\Bigl(e^{\partial_j}\sqrt{x_j(x_k+a_k)}-e^{\partial_k}\sqrt{x_k(x_j+a_j)}\Bigr)\phi_0(\bm{x})\\
&=\sqrt{(x_j+1)(x_k+a_k)}\cdot\sqrt{\frac{(x_j+a_j)(N-|x|)}{(x_j+1)(N-|x|+b-1)}}\,\phi_0(\bm{x})\\
&\ -\sqrt{(x_k+1)(x_j+a_j)}\cdot\sqrt{\frac{(x_k+a_k)(N-|x|)}{(x_k+1)(N-|x|+b-1)}}\,\phi_0(\bm{x})=0.
\end{align*}
The derivation of \eqref{Htherm} goes as follows. Due to \eqref{Hherm}, one gets
\begin{align*}
\langle \mathcal{H}_T\phi_0 p,\phi_0 q\rangle=\langle \phi_0 p,\mathcal{H}_T\phi_0 q\rangle,
\end{align*}
By using the intertwining relation \eqref{intert}, one is led to
\begin{align*}
l.h.s.=\langle\phi_0\widetilde{\mathcal H}_T p,\phi_0 q\rangle=(\widetilde{\mathcal H}_T p,q),\quad
r.h.s.=\langle\phi_0 p,\phi_0\widetilde{\mathcal H}_T q\rangle=(p,\widetilde{\mathcal H}_T q).
\end{align*}
%
%
\section{Eigenpolynomials}
\label{sec:eipoly}
The  complete set of eigenvectors of $\mathcal{H}_T$ \eqref{HTdef} 
consists of products of eigenpolynomials of each $\widetilde{\mathcal H}_i$,  $i=0,1,\ldots,n-1$.
These consist of two types of polynomials, the one depends on  a single variable and the other
depends on two variables.
The first type polynomials are  the eigenpolynomials of $\widetilde{\mathcal H}_0$ and
they are single variable Hahn polynomials.
%
%
\subsection{Single variable Hahn polynomials}
\label{sec:SHahn}
The single variable Hahn polynomial \cite{askey, KLS,os12}
\begin{align}
H_m(x;a,b,N)&\eqdef {}_3F_2\Bigl(
  \genfrac{}{}{0pt}{}{-m,\,m+a+b-1,\,-x}
  {a,\,-N}\Bigm|1\Bigr)
  \label{Hahnpoly}\\[2pt]
  &=\sum_{k=0}^m\frac{(-m)_k(m+a+b-1)_k(-x)_k}{(a)_k(-N)_k\ k!},\nonumber
\end{align}
satisfies a difference equation
\begin{align}
&(N-x)(x+a)\bigl(H_m(x)-H_m(x+1)\bigr)+x(N-x+b)\bigl(H_m(x)-H_m(x-1)\bigr)\n
&\qquad =m(m+a+b-1)H_m(x),
\label{Hahneq}
\end{align}
and they are orthogonal with respect to the hypergeometric binomial distribution
\begin{align*}
&\sum_{x=0}^NH_m(x)H_{m'}(x)W_2(x;\bm{\lambda})=0,\quad m\neq m',\\
&\quad \phi_0(x;\bm{\lambda})^2=W_2(x;\bm{\lambda})=\frac{N!}{x!(N-x)!}\frac{(a)_x(b)_{N-x}}{(a+b)_N},
\quad \bm{\lambda}\eqdef(a,b,N). 
\end{align*}
Here the parameter dependence of $H_m(x)$ is suppressed for simplicity.
The corresponding difference operator is
\begin{align}
\widetilde{\mathcal H}&=B(x)(1-e^\partial)+D(x)(1-e^{-\partial})
=\Bigl(B(x)-D(x)e^{-\partial}\Bigr)\bigl(1-e^{\partial}\bigr),
\label{HaHfac}\\
&\qquad B(x)=(N-x)(x+a),\quad D(x)=x(N-x+b).
\end{align}
A similarity transformation in terms of $\sqrt{W_2(x)}$ produces a real symmetric positive semi-definite
matrix $\mathcal{H}$,
\begin{align}
\mathcal{H}&=\mathcal{A}(\bm{\lambda})^T\mathcal{A}(\bm{\lambda}),
\label{HaAA}\\
&\quad \mathcal{A}(\bm{\lambda})=\sqrt{(N-x)(x+a)}-e^\partial\sqrt{x(N-x+b)},\n
&\quad \mathcal{A}(\bm{\lambda})^T=\sqrt{(N-x)(x+a)}-\sqrt{x(N-x+b)}\,e^{-\partial}.\nonumber
\end{align}
It is well known that this system is {\em shape invariant} \cite{os12}(4.2),
\begin{equation}
  \mathcal{A}(\bm{\lambda})\mathcal{A}(\bm{\lambda})^T
  =\mathcal{A}(\bm{\lambda}+\bm{\delta})^T
  \mathcal{A}(\bm{\lambda}+\bm{\delta})+\mathcal{E}(1;\bm{\lambda}),\quad
  \mathcal{E}(1;\bm{\lambda})=a+b,\quad \bm{\delta}=(1,1,-1).
  \label{shapeinv1}
\end{equation}
Two important consequences of the shape invariance are the eigenvalue formula \cite{os12}(4.5)
\begin{align}
  &\mathcal{E}(m\,;\bm{\lambda})
  =\sum_{s=0}^{m-1}\mathcal{E}(1\,;\bm{\lambda}+s\bm{\delta})=\sum_{s=0}^{m-1}(a+s+b+s)=m(m+a+b-1).
  \label{spectrumform}
\end{align}
and the Rodrigues type formula \cite{os12}(4.6)
\begin{align}
  &H_m(x\,;\bm{\lambda})\phi_0(x;\bm{\lambda})\propto
  \mathcal{A}(\bm{\lambda})^T
  \mathcal{A}(\bm{\lambda}+\bm{\delta})^T\cdots
  \mathcal{A}\bigl(\bm{\lambda}+(m-1)\bm{\delta}\bigr)^T
  \phi_0(x\,;\bm{\lambda}+m\bm{\delta}).
  \label{eigvecform}
\end{align}
These two formulas \eqref{spectrumform} and \eqref{eigvecform} are essential for the determination of 
the second type eigenpolynials and the eigenvalues to be derived in \S\ref{sec:hieig}.
The second type polynomials are the eigenpolynomials of $\widetilde{\mathcal H}_i$, $i=1,\ldots,n-1$
and they are closely related to the `renormalised Hahn polynomials' \cite{dunkl0}(3.1).
Corresponding to the factorisation of $\widetilde{H}$ \eqref{HaHfac}, $H_m(x)$ satisfies the forward \cite{KLS}(9.5.7) \cite{os12}(4.20) and backward 
\cite{KLS}(9.5.8) \cite{os12}(4.21) shift relations,
\begin{align}
\text{Forward shift:}\quad &H_m(x,a,b,N)-H_m(x-1,a,b,N)\n
&\qquad =\frac{m(m+a+b-1)}{aN}H_{m-1}(x,a+1,b+1,N-1),
\label{Hfor}\\
\text{Backward shift:}\quad  &(N-x)(x+a)H_m(x,a+1,b+1,N-1)\n
& -x(N-x+b)H_m(x-1,a+1,b+1,N-1)\n
& =aNH_{m+1}(x,a,b,N).
\label{Hback}
\end{align}

%
%
\subsection{Type one eigenpolynomials}
\label{sec:typeone}

Let us apply the difference operator $\widetilde{\mathcal H}_T$ 
to an $m$ degree polynomial in $x_i$, $p_m(x_i)$,
\begin{align*}
\widetilde{\mathcal H}_Tp_m(x_i)
&=\!(N-|x|)(x_i+a_i)\bigl(p_m(x_i)-p_m(x_i+1)\bigr)\!+\!x_i(N-|x|+b)\bigl(p_m(x_i)-p_m(x_i-1)\bigr)\\
&\ +\sum_{j\neq i}x_j(x_i+a_i)\bigl(p_m(x_i)-p_m(x_i+1)\bigr)
   +x_i\sum_{j\neq i}(x_j+a_j)\bigl(p_m(x_i)-p_m(x_i-1)\bigr)\\
&=(N-x_i)(x_i+a_i)\bigl(p_m(x_i)-p_m(x_i+1)\bigr)\\
&\quad +x_i(N-x_i+|a|+b-a_i)\bigl(p_m(x_i)-p_m(x_i-1)\bigr).
\end{align*}
This means that $H_m(x_i;a_i,|a|+b-a_i,N)$ is an eogenpolynomial of $\widetilde{\mathcal H}_T$
with the eigenvalue $\mathcal{E}_T(m)=m(m+|a|+b-1)$.
Likewise, let us apply the difference operator $\widetilde{\mathcal H}_T$ 
to an $m$ degree polynomial in $|x|=x_1+\ldots+x_n$, $p_m(|x|)$. In this case, 
as $e^{-\partial_j+\partial_k}p_m(|x|)=p_m(|x|)$, $\widetilde{\mathcal H}_1$ \eqref{Ht1def}
annihilates $p_m(|x|)$, $\widetilde{\mathcal H}_1p_m(|x|)=0$, 
and only $\widetilde{\mathcal H}_0$ \eqref{Ht0def} acts on $p_m(|x|)$,
\begin{align*}
\widetilde{\mathcal H}_Tp_m(|x|)&=\widetilde{\mathcal H}_0p_m(|x|)\\
&=(N-|x|)(|x|+|a|)\bigl(p_m(|x|)-p_m(|x|+1)\bigr)\\
&\quad +|x|(N-|x|+b)\bigl(p_m(|x|)-p_m(|x|-1)\bigr).
\end{align*}
This means that $H_m(|x|;|a|,b,N)$ is an eigenpolynomial of $\widetilde{\mathcal H}_T$
with the same eigenvalue $\mathcal{E}_T(m)=m(m+|a|+b-1)$ as above.
Let us introduce  a variable $x_J$ and a parameter $a_J$ by
\begin{equation}
\text{Definition}: \quad x_J\eqdef\sum_{j\in J}x_j,\quad a_J\eqdef \sum_{j\in J}a_j,
\label{xjdef}
\end{equation}
in which the set $J\subseteq\{1,2,\ldots,n\}$ consists of {\em distinct elements\/}.
The same calculation with $\widetilde{\mathcal H}_Tp_m(x_J)$ leads to the following
\begin{prop}
\label{prop:single}
There are many {\rm (at least }$2^n-1${\rm)}  eigenpolynomials of degree $m$ of
$\widetilde{\mathcal H}_T$ with the same eigenvalue 
\begin{align}
\widetilde{\mathcal H}_TH_{m}(x_J;a_J,|a|+b-a_J,N)=m(m+|a|+b-1)H_{m}(x_J;a_J,|a|+b-a_J,N).
\label{HxJ}
\end{align}
With different degrees, they are orthogonal with each other, but the same degree ones are not orthogonal, in general.
\end{prop}
\begin{rema}
\label{rema:sameeig}
These results would suggest that all  the total degree $M$ eienpolynomials constituting 
the complete set of eigenvectors of $\mathcal{H}_T$ have the same eigenvalue $\mathcal{E}_T(M)=M(M+|a|+b-1)$.
This will be shown in  {\bf Main Theorem \ref{theo:main}}. 
It should be stressed that this is a very special situation 
and it is due to the exceptionally symmetric difference operator $\widetilde{\mathcal H}_T$.
The situation is very different from other discrete multivariate orthogonal polynomials satisfying
difference equations {\rm \cite{I23,mKrawt, mRah}}, in which  the degree one eigenvalues and eigenpolynomials are irrational.

\end{rema}
%
%
\subsection{Type two eigenpolynomials}
\label{sec:typetwo}
Let us introduce a family of variables $x_{>i}$   and the corresponding parameters $a_{>i}$.
\begin{equation}
\text{Definition:}\quad  x_{>i}\eqdef x_{i+1}+\cdots+x_n,\quad
a_{>i}\eqdef a_{i+1}+\cdots+a_n,\quad i=1.\ldots,n-1.
\label{x>def}
\end{equation}
By simple calculations one arrives at the following
\begin{prop}
\label{prop:tdef}
Degree one type two eigenpolynomials $\{t_i\}$, $i=1,\ldots,n-1$ are orthogonal with each other.
\begin{align}
&t_i\eqdef\widetilde{\mathcal H}_ix_i=a_{>i}x_i-a_ix_{>i},\quad
\widetilde{\mathcal H}_it_i=(a_i+a_{>i})t_i,\quad i=1,\ldots,n-1,
\label{tidef}\\
& \widetilde{\mathcal H}_T\,t_i=(|a|+b)t_i,\qquad \widetilde{\mathcal H}_j\,t_i=0,\quad j>i,
\quad \Rightarrow (t_i,t_j)=0, \quad i\neq j.
\label{ti2}
\end{align}
\end{prop}
Since a degree one type one eigenpolynomial of $\widetilde{\mathcal H}_T$,  
$h_1(x_i)\eqdef H_1(x_i;a_i,|a|+b-a_i,N)$ has a form 
$h_1(x_i)=1-(|a|+b)x_i/(a_iN)$, and $\widetilde{\mathcal H}_i1=0$, 
$t_i$ can be expressed as $t_i=\alpha\widetilde{\mathcal H}_ih_1(x_i)$, $\alpha=-a_iN/(|a|+b)$.
Therefore $t_i$ is also an eigenpolynomial of $\widetilde{\mathcal H}_T$ as
\begin{equation*}
\widetilde{\mathcal H}_Tt_i=\alpha \widetilde{\mathcal H}_T\widetilde{\mathcal H}_ih_1(x_i)
=\alpha\widetilde{\mathcal H}_i\widetilde{\mathcal H}_Th_1(x_i)
=(|a|+b)\alpha\widetilde{\mathcal H}_ih_1(x_i)=(|a|+b)t_i.
\end{equation*}
The annihilation of $t_i$ by $\widetilde{\mathcal H}_j$ ($j>i$) 
is the consequence that $\widetilde{\mathcal H}_j$ ($j>i$)
contains variables $x_{j},x_{j+1}\ldots,x_n$. Therefore $1-e^{-\partial_j+\partial_k}$ 
annihilates $x_i$ and $x_{>i}$. 
The orthogonality is the consequence of \eqref{Htherm}
\begin{equation*}
(\widetilde{\mathcal H}_jt_j,t_i)=(a_j+a_{>j})(t_j,t_i)=(t_j,\widetilde{\mathcal H}_jt_i)=0.
\end{equation*}

It is expected that the higher degree eigenpolynomials of $\widetilde{\mathcal H}_i$ would be
orthogonal with the counterparts of $\widetilde{\mathcal H}_j$. 
The next task is to construct the eigenpolynomials of $\widetilde{\mathcal H}_i$ depending on two
variables $x_i$ and $x_{>i}$, $i=1,\ldots,n-1$.
%
%
\subsection{Eigenpolynomials of $\widetilde{\mathcal H}_i$ }
\label{sec:hieig}
The difference operator $\widetilde{\mathcal H}_i$,
\begin{equation*}
\widetilde{\mathcal H}_i=\sum_{j\neq k\ge i}x_j(x_k+a_k)(1-e^{-\partial_j+\partial_k}),
\tag{\ref{Htidef}}
\end{equation*}
is greatly simplified when its action is restricted to a polynomial space $V_N^{(i)}(\bm{x})$ spanned by
two variables $x_i$ and $x_{>i}$
\begin{align}
V_N^{(i)}(\bm{x})&
={\rm Span}\left\{x_i^{m_i}x_{>i}^{m_{>i}}\mid 0\le m_i+m_{>i} \le N\right\}\subseteq V_N(\bm{x}),
\label{Vres}\\
\widetilde{\mathcal H}_i\to \widetilde{\mathcal H}^{(i)}&\eqdef
x_i(x_{>i}+a_{>i})\bigl(1-e^{-\partial_i+\partial_{>i}}\bigr)
+x_{>i}(x_i+a_i)\bigl(1-e^{\partial_i-\partial_{>i}}\bigr),
\label{resHti}\\
&=\Bigl(x_{>i}(x_i+a_i)e^{-\partial_{>i}}-x_i(x_{>i}+a_{>i})e^{-\partial_i}\Bigr)\bigl(e^{\partial_{>i}}-e^{\partial_i}\bigr).
\label{Hifact}
\end{align}
The corresponding real symmetric matrix $\mathcal{H}^{(i)}$ is obtained by a similarity transformation
\eqref{intert} in terms of $\phi_0^{(i)}$, the square root of the corresponding weight function, 
\begin{align}
&{\phi_0(\bm{x};\bm{\lambda}) ^{(i)2}}=W^{(i)}(\bm{x};\bm{\lambda}) 
\qquad \bm{\lambda}=(a_i,a_{>i},b',N),\quad  b'\eqdef b+|a|-a_i-a_{>i},\n
&\qquad \qquad \ \  =\frac{N!}{x_i!x_{>i}!(N-x_i-x_{>i})!}\frac{(a_i)_{x_i}(a_{>i})_{x_{>i}}(b')_{N-x_i-x_{>i}}}{(|a|+b)_N},
\label{Wi}\\
\mathcal{H}^{(i)}&\eqdef\phi_0^{(i)}(\bm{x};\bm{\lambda})
\widetilde{\mathcal H}^{(i)}\phi_0^{(i)}(\bm{x};\bm{\lambda})^{-1},\n
&=x_i(x_{>i}+a_{>i})+x_{>i}(x_i+a_i)-\sqrt{x_i(x_i+a_i-1)(x_{>i}+1)(x_{>i}+a_{>i})}\,e^{-\partial_i+\partial_{>i}}\n
&\hspace{50mm}
-\sqrt{(x_i+1)(x_i+a_i)x_{>i}(x_{>i}+a_{>i}-1)}\,e^{\partial_i-\partial_{>i}}\\
&=\mathcal{B}({\bm \lambda})^{(i)T}\mathcal{B}({\bm \lambda})^{(i)}
\label{Hupi1}\\
&\quad \mathcal{B}({\bm \lambda})^{(i)}\eqdef 
e^{\partial_{>i}}\sqrt{x_{>i}(x_i+a_i)}-e^{\partial_i}\sqrt{x_i(x_{>i}+a_{>i})},\n
&\quad \mathcal{B}({\bm \lambda})^{(i)T}
= \sqrt{x_{>i}(x_i+a_i)}\,e^{-\partial_{>i}}-\sqrt{x_i(x_{>i}+a_{>i})}\,e^{-\partial_i},\n
&\mathcal{B}({\bm \lambda})^{(i)}\phi_0^{(i)}(\bm{x};\bm{\lambda})=0\quad \Rightarrow
\mathcal{H}^{(i)}\phi_0^{(i)}(\bm{x};\bm{\lambda})=0.\nonumber
\end{align}
The similar structure of $\mathcal{H}^{(i)}$ \eqref{Hupi1} with $\mathcal{H}$ of 
the single variable Hahn \eqref{HaAA} leads to the
following
\begin{theo}
\label{theo:Pishape}{\bf Shape invariance} determines the eigenpolynomial 
$\mathcal{P}_m^{(i)}(\bm{x};\bm{\lambda})$ of $\widetilde{\mathcal H}^{(i)}$
and its eigenvalue $\mathcal{E}^{(i)}(m)$. They are orthogonal with each other.
\begin{align}
 \mathcal{P}_{m}^{(i)}\bigl(\bm{x};\bm{\lambda})&= \mathcal{P}_{m}^{(i)}\bigl(x_i,x_{>i};a_i,a_{>i}\bigr)\n
\qquad\qquad &\eqdef
 \sum_{k=0}^{m}(-1)^k\binom{m}{k}(a_{>i}+k)_{m-k}(a_i+m-k)_k(-x_i)_{m-k}(-x_{>i})_k,
 \label{Piform}\\
\widetilde{\mathcal H}_i \mathcal{P}_{m}^{(i)}\bigl(\bm{x};\bm{\lambda}\bigr)
&=\widetilde{\mathcal H}^{(i)} \mathcal{P}_{m}^{(i)}\bigl(\bm{x};\bm{\lambda}\bigr)
=\mathcal{E}_{i}(m)\mathcal{P}_{m}^{(i)}\bigl(\bm{x};\bm{\lambda}\bigr),\n
\qquad \mathcal{E}_{i}(m)&=m(m+a_i+a_{>i}-1),
\label{eigHi}\\
\widetilde{\mathcal H}_j \mathcal{P}_{m}^{(i)}\bigl(\bm{x};\bm{\lambda}\bigr)&=0,\quad j>i
\quad \Rightarrow (\mathcal{P}_{m}^{(i)},\mathcal{P}_{m}^{(j)})=0,\quad i\neq j.
\label{Pmortho}
\end{align}
\end{theo}
It is straightforward to verify the shape invariance. 
Derivation of the eigenvalue formula is the same as that of the single variable Hahn \eqref{spectrumform}.
The Rodrigues like formula reads as below.
\begin{align}
& \mathcal{B}(\bm{\lambda})\mathcal{B}(\bm{\lambda})^T
  =\mathcal{B}(\bm{\lambda}+\bm{\delta})^T
  \mathcal{B}(\bm{\lambda}+\bm{\delta})+a_i+a_{>i},\quad
 \quad \bm{\delta}=(1,1,0,-1),
  \label{shapeinvB}\\
  & \Rightarrow \mathcal{P}_m^{(i)}(\bm{x};\bm{\lambda})
  \phi_0^{(i)}(\bm{x};\bm{\lambda})\propto\mathcal{B}(\bm{\lambda})^T
  \mathcal{B}(\bm{\lambda}+\bm{\delta})^T\cdots
  \mathcal{B}\bigl(\bm{\lambda}+(m-1)\bm{\delta}\bigr)^T
  \phi_0(x\,;\bm{\lambda}+m\bm{\delta}),
  \label{type2Rod}\\
 & \mathcal{P}_0^{(i)}=1,\quad  \mathcal{P}_1^{(i)}(\bm{x};\bm{\lambda})=a_ix_{>i}-a_{>i}x_i=-t_i.\nonumber
\end{align}
 Appendix A provides an inductive proof of the explicit expression of the type two eigenpolynomials
 \eqref{Piform} based on the Rodrigues type formula \eqref{type2Rod}.
 \begin{rema}
 \label{rema:Genest}
 The type two eigenpolynomial \eqref{Piform} for the special case of $n=2$ theory appeared in
 Genest-Vinet paper \rm{\cite{genest}}. In my notation it was up to an overall sign $(-1)^m$
 \begin{equation}
(a_i)_m(-x_i-x_{>i})_mH_m(x_i;a_i,a_{>i},-x_i-x_{>i}).
\label{anotherHahn}
\end{equation}
This definition of Hahn polynomial was also used by many, for example, Tratnik \cite{tra3}.
I arrived at \eqref{Piform} by constructing the eigenpolynomials  of $n=2$ theory based on 
$\widetilde{\mathcal H}_0$ \eqref{Ht0def} by hand up to degree 4.
 \end{rema}
 \begin{rema}
 \label{rema:special}
 The type two eigenpolynomial \eqref{Piform} has a simple special value at $x_i=m$, $x_{>i}=0$,
 \begin{equation}
 \mathcal{P}_{m}^{(i)}\bigl(m,0;a_i,a_{>i}\bigr)=(-1)^mm!(a_{>i})_m.
 \label{sPval}
 \end{equation}
 \end{rema}
 \begin{prop}
 \label{prop:PiForbac}
 In parallel with the single variable case \eqref{Hfor}, \eqref{Hback}, 
 forward and backward shift relations hold for the type two
 eigenpolynomials \eqref{Piform} corresponding to the factorisation of 
 $\widetilde{\mathcal H}^{(i)}$ \eqref{Hifact}, 
 \begin{align}
\text{Forward shift:}\quad &\mathcal{P}_{m}^{(i)}\bigl(x_i,x_{>i}+1;a_i,a_{>i}\bigr)
-\mathcal{P}_{m}^{(i)}\bigl(x_i+1,x_{>i};a_i,a_{>i}\bigr)\n
&\quad  =m(m+a_i+a_{>i}-1)\mathcal{P}_{m-1}^{(i)}\bigl(x_i,x_{>i};a_i+1,a_{>i}+1\bigr),
 \label{Pifor}\\
\text{Backward shift:}\quad   &x_{>i}(x_i+a_i)\mathcal{P}_{m}^{(i)}\bigl(x_i,x_{>i}-1;a_i+1,a_{>i}+1\bigr)\n
&\qquad -x_i(x_{>i}+a_{>i})\mathcal{P}_{m}^{(i)}\bigl(x_i-1,x_{>i};a_i+1,a_{>i}+1\bigr)\n
&\quad =\mathcal{P}_{m+1}^{(i)}\bigl(x_i,x_{>i};a_i,a_{>i}\bigr).
\label{Piback}
 \end{align}
 \end{prop}
 They are the consequences of the shape invariance \eqref{shapeinvB} and the factorisation of 
 $\widetilde{\mathcal H}^{(i)}$ \eqref{Hifact}. For the forward shift, it is known \cite{os12}(4.20)
 \begin{equation*}
\bigl(e^{\partial_{>i}}-e^{\partial_i}\bigr)\mathcal{P}_{m}^{(i)}\bigl(x_i,x_{>i};a_i,a_{>i}\bigr)
\propto \mathcal{P}_{m-1}^{(i)}\bigl(x_i,x_{>i};a_i+1,a_{>i}+1\bigr).
\end{equation*}
Evaluation of  both sides at $x_i=m$, $x_{>i}=0$ determines the coefficient $m(m+a_i+a_{>i}-1)$, the eigenvalue.
Similarly for the backward shift relation,
\begin{align*}
&\Bigl(x_{>i}(x_i+a_i)e^{-\partial_{>i}}-x_i(x_{>i}+a_{>i})e^{-\partial_i}\Bigr)
\mathcal{P}_{m}^{(i)}\bigl(x_i,x_{>i};a_i+1,a_{>i}+1\bigr)\\
&\qquad \propto \mathcal{P}_{m+1}^{(i)}\bigl(x_i,x_{>i};a_i,a_{>i}\bigr)
\end{align*}
is known \cite{os12}(4.21) . Evaluation at $x_i=m+1$, $x_{>i}=0$ sets the equality.
 Two types of forward and backward shift relations are reported in \cite{genest}A.6,A.7 for the total 
 eigenpolynomials. The one is for the type two and the other is for the single variable Hahn.
 \begin{theo}
 \label{theo:rec}
 The type two eigenpolynomials \eqref{Piform} satisfy the following recursion relations.
 \begin{align}
 \text{Forward recursion:}\quad &(x_{i}+a_{i})\mathcal{P}_{m}^{(i)}\bigl(x_i+1,x_{>i};a_i,a_{>i}\bigr)
                                                           +(x_{>i}+a_{>i})\mathcal{P}_{m}^{(i)}\bigl(x_i,x_{>i}+1;a_i,a_{>i}\bigr)\n                                                        
 &\quad =(x_i+x_{>i}+a_i+a_{>i}+m) \mathcal{P}_{m}^{(i)}\bigl(x_i,x_{>i};a_i,a_{>i}\bigr),
\label{rec1}\\
 \text{Backward recursion:}\quad & x_{i}\mathcal{P}_{m}^{(i)}\bigl(x_i-1,x_{>i};a_i,a_{>i}\bigr)
                                            +x_{>i}\mathcal{P}_{m}^{(i)}\bigl(x_i,x_{>i}-1;a_i,a_{>i}\bigr)\n  
&\quad =(x_i+x_{>i}-m) \mathcal{P}_{m}^{(i)}\bigl(x_i,x_{>i};a_i,a_{>i}\bigr).
\label{rec2}
\end{align}                                                         
 \end{theo}
 For the proof of the Theorem, let us define two difference operators
 \begin{align}
\text{Definition:}\quad \mathcal{F}^{(i)}\eqdef (x_i+a_i)e^{\partial_i}+(x_{>i}+a_{>i})e^{\partial_{>i}},
\quad  \mathcal{G}^{(i)}\eqdef x_i\,e^{-\partial_i}+x_{>i}\,e^{-\partial_{>i}},
 \label{FGdef}
\end{align}
Apply $\mathcal{G}^{(i)}$ to the variable dependent part of the $k$-th term of $\mathcal{P}_{m}^{(i)}\bigl(x_i,x_{>i};a_i,a_{>i}\bigr)$ \eqref{Piform},
 \begin{align}
 \mathcal{G}^{(i)}\,(-x_i)_{m-k}(-x_{>i})_k&=x_i(-x_i+1)_{m-k}(-x_{>i})_k+(-x_i)_{m-k}x_{>i}(-x_{>i}+1)_k\n
 &=(x_i+x_{>i}-m)\,(-x_i)_{m-k}(-x_{>i})_k,
 \end{align}
 and the backward recursion \eqref{rec2} is proved. It is interesting to note that 
 $\widetilde{\mathcal H}^{(i)}+\mathcal{F}^{(i)}\mathcal{G}^{(i)}$ is a function only,
 \begin{align*}
&\widetilde{\mathcal H}^{(i)}+\mathcal{F}^{(i)}\mathcal{G}^{(i)}=(x_i+x_{>i}+1)(x_i+x_{>i}+a_i+a_{>i}).\\
&\Rightarrow \mathcal{F}^{(i)}\mathcal{G}^{(i)}\mathcal{P}_{m}^{(i)}=\Bigl((x_i+x_{>i}+1)(x_i+x_{>i}+a_i+a_{>i})-m(m+a_i+a_{>i}-1)\Bigr)\mathcal{P}_{m}^{(i)}. \tag{*}
\end{align*}
By using \eqref{rec2}, $\mathcal{F}^{(i)}\mathcal{G}^{(i)}\mathcal{P}_{m}^{(i)}$ is evaluated as
\begin{equation*}
\mathcal{F}^{(i)}\mathcal{G}^{(i)}\mathcal{P}_{m}^{(i)}=(x_i+x_{>i}+1-m)\mathcal{F}^{(i)}\mathcal{P}_{m}^{(i)}.
\tag{**}
\end{equation*}
Combining $(*)$ and $(**)$ proves \eqref{rec1}.
The forward and backward recursions (Theorem \ref{theo:rec})  lead to the following
\begin{theo}
\label{theo:eigens}
{\bf Eigenvalue Theorem}
\begin{align}
\widetilde{\mathcal H}_T\mathcal{P}_{m}^{(i)}\bigl(x_i,x_{>i};a_i,a_{>i}\bigr)&
=m(m+|a|+b-1)\mathcal{P}_{m}^{(i)}\bigl(x_i,x_{>i};a_i,a_{>i}\bigr),
\label{eigT}\\
\widetilde{\mathcal H}_j\mathcal{P}_{m}^{(i)}\bigl(x_i,x_{>i};a_i,a_{>i}\bigr)&
=m(m+a_j+\cdots+a_n-1)\mathcal{P}_{m}^{(i)}\bigl(x_i,x_{>i};a_i,a_{>i}\bigr),\quad j<i,
\label{eigj}\\
&=m(m+a_j+a_{>j}-1)\mathcal{P}_{m}^{(i)}\bigl(x_i,x_{>i};a_i,a_{>i}\bigr),\qquad \quad j<i,\n
\widetilde{\mathcal H}_i\mathcal{P}_{m}^{(i)}\bigl(x_i,x_{>i};a_i,a_{>i}\bigr)&
=m(m+a_i+\cdots+a_n-1)\mathcal{P}_{m}^{(i)}\bigl(x_i,x_{>i};a_i,a_{>i}\bigr),
\label{eigi}\\
&=m(m+a_i+a_{>i}-1)\mathcal{P}_{m}^{(i)}\bigl(x_i,x_{>i};a_i,a_{>i}\bigr),\n
\widetilde{\mathcal H}_j\mathcal{P}_{m}^{(i)}\bigl(x_i,x_{>i};a_i,a_{>i}\bigr)&=0,\quad j>i,
\hspace{50mm}
i=1,\ldots,n-1.
\label{eig0}
\end{align}
The first (m=1) eigenvalue is always the sum of the parameters $a_i$ and $b$ contained in the difference operator
$\widetilde{\mathcal H}$.
\end{theo}
Simple proofs of eigenvalue formulas of $\widetilde{\mathcal H}_j$ \eqref{eigj} 
and $\widetilde{\mathcal H}_T$ \eqref{eigT} are relegated in Appendix B.
\section{General eigenpolynomials}
\label{sec:geneig}
Let us start with 
\begin{theo}
\label{theo:main}{\bf Main Theorem}
The complete set of eigenpolynomials of the multivariate Hahn system defined by {\bf Definition \ref{def:Ht0}}
is
\begin{align}
\mathcal{P}_{\bm m}(\bm{x};\bm{a},b,N)
&=\prod_{i=1}^{n-1}\mathcal{P}_{m_i}^{(i)}
\Bigl(x_i,x_{>i}-\sum_{j=i+1}^{n-1}m_j;a_i,a_{>i}+2\sum_{j=i+1}^{n-1}m_j\Bigr)\n
& \quad \times H_{m_0}\Bigl(|x|-\sum_{j=1}^{n-1}m_j;|a|+2\sum_{j=1}^{n-1}m_j,b,N-\sum_{j=1}^{n-1}m_j\Bigr)
\label{sol1}\\
&=\prod_{i=1}^{n-1}\left( \sum_{k=0}^{m_i}(-1)^k\binom{m_i}{k}
\Bigl(a_{>i}+2\sum_{j=i+1}^{n-1}m_j+k\Bigr)_{m_i-k}
\Bigl(a_i+m_i-k\Bigr)_k\right.\n
&\hspace{20mm} \times\left. \Bigl(-x_i\Bigr)_{m_i-k}\Bigl(-x_{>i}+\sum_{j=i+1}^{n-1}m_j\Bigr)_k\right)\n
&\quad \times {}_3F_2\left(
  \genfrac{}{}{0pt}{}{-m_0,\,m_0+2\sum_{j=1}^{n-1}m_j+|a|+b-1,\,-|x|+\sum_{j=1}^{n-1}m_j}
  {|a|+2\sum_{j=1}^{n-1}m_j,\,-N+\sum_{j=1}^{n-1}m_j}\Bigm|1\right)
  \label{sol2}\\[4pt]
& \hspace{-20mm}\widetilde{\mathcal H}_T\mathcal{P}_{\bm m}(\bm{x};\bm{a},b,N)
=|m|(|m|+|a|+b-1)\mathcal{P}_{\bm m}(\bm{x};\bm{a},b,N),
\label{Hteq}\\
& \hspace{-20mm}\widetilde{\mathcal H}_i\mathcal{P}_{\bm m}(\bm{x};\bm{a},b,N)
=|m|\Bigl(|m|+\sum_{j=i}^{n}a_j-1\Bigr)\mathcal{P}_{\bm m}(\bm{x};\bm{a},b,N),\quad i=1,\ldots, n-1,
\label{Htieq}\\
&\hspace{20mm}  \bm{m}=(m_0,m_1,\ldots,m_{n-1})\in\cX,\n
&\Big(\mathcal{P}_{\bm m},\mathcal{P}_{\bm{m}'}\Bigr)=0,\qquad \bm{m}\neq\bm{m}'.
\label{genortho}
\end{align}
The numbering of $\bm{m}$ is different from \eqref{VNdef} but no confusion will arise.
\end{theo}
The proof in several steps follows. The first task is to combine two adjacent type two eigenpolynomials
$\mathcal{P}_{m_i}^{(i)}$ and $\mathcal{P}_{m_{i-1}}^{(i-1)}$  to form a product
eigenpolynomial of $\widetilde{\mathcal H}_{i-1}$ of degree $m_i+m_{i-1}$.
The process starts from the top one $\mathcal{P}^{(n-1)}$.  It is achieved by adjusting the 
variable and the parameter of
the lower one in terms of the degree of the above one as stated in the following
\begin{theo}
\label{theo:ii-1}
\begin{align}
&\widetilde{\mathcal H}_{i-1}\mathcal{P}_{m_i}^{(i)}\bigl(x_i,x_{>i};a_i,a_{>i}\bigr)\mathcal{P}_{m_{i-1}}^{(i-1)}\bigl(x_{i-1},x_{>i-1}-m_i;a_i,a_{>i-1}+2m_i\bigr)\n
&\quad =\bigl(m_i+m_{i-1}\bigr)\bigl(m_i+m_{i-1}+a_{i-1}+\ldots+a_n-1)\n
&\qquad\qquad \times\mathcal{P}_{m_i}^{(i)}\bigl(x_i,x_{>i};a_i,a_{>i}\bigr)\mathcal{P}_{m_{i-1}}^{(i-1)}\bigl(x_{i-1},x_{>i-1}-m_i;a_{i-1},a_{>i-1}+2m_i\bigr).
\label{sumii-1}
\end{align}
The combined eigenpolynomial satisfies the Eigenvalue Theorem \ref{theo:eigens}.
\end{theo}
The proof is slightly lengthy but I show it here as it is an essential part of this paper,
\begin{align*}
\widetilde{\mathcal H}_{i-1}=&\widetilde{\mathcal H}_{i}+\widetilde{\mathcal H}_{idif},\n
\widetilde{\mathcal H}_{idif}=&\ x_{i-1}(x_i+a_i)(1-e^{-\partial_{i-1}+\partial_i})\n
& +x_{i-1}(x_{>i}+a_{>i})(1-e^{-\partial_{i-1}+\partial_{>i}})\n
& +x_i(x_{i-1}+a_{i-1})(1-e^{\partial_{i-1}-\partial_i})\n
& +x_{>i}(x_{i-1}+a_{i-1})(1-e^{\partial_{i-1}-\partial_{>i}}).
\end{align*}
Let us apply $\widetilde{\mathcal H}_{idif}$ to $\mathcal{P}_{m_i}^{(i)}\mathcal{P}_{m_{i-1}}^{(i-1)}$ with the 
above arguments. The default arguments are suppressed and the changed parts only displayed. 
Since $\mathcal{P}_{m_i}^{(i)}$ does not contain $x_{i-1}$, it reads
\begin{align*}
&\widetilde{\mathcal H}_{idif}\mathcal{P}^{(i)}\mathcal{P}^{(i-1)}\\
&\quad =\quad x_{i-1}(x_i+a_i)\Bigl(\mathcal{P}^{(i)}\mathcal{P}^{(i-1)}-\mathcal{P}^{(i)}(x_i+1)\mathcal{P}^{(i-1)}(x_{i-1}-1,x_{>i-1}+1)\Bigr)\\
&\quad +x_{i-1}(x_{>i}+a_{>i})\Bigl(\mathcal{P}^{(i)}\mathcal{P}^{(i-1)}-\mathcal{P}^{(i)}(x_{>i}+1)\mathcal{P}^{(i-1)}(x_{i-1}-1,x_{>i-1}+1)\Bigr)\\
&\quad +x_i(x_{i-1}+a_{i-1})\Bigl(\mathcal{P}^{(i)}\mathcal{P}^{(i-1)}-\mathcal{P}^{(i)}(x_i-1)\mathcal{P}^{(i-1)}(x_{i-1}+1,x_{>i-1}-1)\Bigr)\\
&\quad +x_{>i}(x_{i-1}+a_{i-1})\Bigl(\mathcal{P}^{(i)}\mathcal{P}^{(i-1)}-\mathcal{P}^{(i)}(x_{>i}-1)\mathcal{P}^{(i-1)}(x_{i-1}+1,x_{>i-1}-1)\Bigr).
\end{align*}
In the first two lines, the changed $\mathcal{P}^{(i-1)}$ are the same. The forward recursion \eqref{rec1} for 
$\mathcal{P}^{(i)}$ can be applied. In the third and fourth lines, the changed $\mathcal{P}^{(i-1)}$ are the same. 
The backward recursion \eqref{rec2} for 
$\mathcal{P}^{(i)}$ can be applied. The result reads
\begin{align*}
&\widetilde{\mathcal H}_{idif}\mathcal{P}^{(i)}\mathcal{P}^{(i-1)}\\
&\quad =x_{i-1}(x_{>i-1}+a_{>i-1})\mathcal{P}^{(i)}\mathcal{P}^{(i-1)} \tag{1}\\
&\quad -x_{i-1}(x_{>i-1}+a_{>i-1}+m_i)\mathcal{P}^{(i)}\mathcal{P}^{(i-1)}(x_{i-1}-1,x_{>i-1}+1) \tag{2}\\
&\quad +x_{>i-1}(x_{i-1}+a_{i-1})\mathcal{P}^{(i)}\mathcal{P}^{(i-1)} \tag{3}\\
&\quad -(x_{>i-1}-m_i)(x_{i-1}+a_{i-1})\mathcal{P}^{(i)}\mathcal{P}^{(i-1)}(x_{i-1}+1,x_{>i-1}-1). \tag{4}
\end{align*}
$(3)+(4)$ is changed to 
\begin{align*}
&(x_{>i-1}-m_i)(x_{i-1}+a_{i-1})\Bigl(\mathcal{P}^{(i)}\mathcal{P}^{(i-1)}-\mathcal{P}^{(i)}\mathcal{P}^{(i-1)}(x_{i-1}+1,x_{>i-1}-1)\Bigr)\n
&\qquad +m_i(x_{i-1}+a_{i-1})\mathcal{P}^{(i)}\mathcal{P}^{(i-1)}.
\end{align*}
Adding $m_ix_{i-1}\mathcal{P}^{(i)}\mathcal{P}^{(i-1)}$ to $(1)$ leads to
\begin{align*}
&\widetilde{\mathcal H}_{idif}\mathcal{P}^{(i)}\mathcal{P}^{(i-1)}\\
&\quad =\mathcal{P}^{(i)}\Bigl[x_{i-1}(x_{>i-1}-m_i +2m_i+a_{>i-1})\Bigl(\mathcal{P}^{(i-1)} -\mathcal{P}^{(i-1)}(x_{i-1}-1,x_{>i-1}+1)\Bigr)\\
&\qquad \qquad \ +(x_{>i-1}-m_i)(x_{i-1}+a_{i-1})
\Bigl(\mathcal{P}^{(i-1)}-\mathcal{P}^{(i-1)}(x_{i-1}+1,x_{>i-1}-1)\Bigr)\Bigr]\\
&\quad +m_ia_{i-1}\mathcal{P}^{(i)}\mathcal{P}^{(i-1)}.
\end{align*}
The difference operator acting on the $\mathcal{P}^{(i-1)}$ is $\widetilde{\mathcal H}_{i-1}$ for the variable
$x_{>i-1}-m_i$ and the argument $a_{>i-1}+2m_i$, leading to
\begin{align*}
&\widetilde{\mathcal H}_{idif}\mathcal{P}^{(i)}\mathcal{P}^{(i-1)}=\Bigl(m_{i-1}(m_{i-1}+a_{i-1}+a_{>i-1}+2m_i-1)+m_ia_{i-1}\Bigr)\mathcal{P}^{(i)}\mathcal{P}^{(i-1)}.
\end{align*}
Adding the contribution from $\widetilde{\mathcal H}_i\mathcal{P}^{(i)}\mathcal{P}^{(i-1)}$ produces
the result \eqref{sumii-1}.\hfill $\square$

Starting from $\mathcal{P}^{(n-1)}$, the combination goes as far as  $\mathcal{P}^{(i)}$ to produce
\begin{align}
\mathcal{R}^{(i)}(\bm{x};\bm{a})&=\prod_{j=i}^{n-1}\mathcal{P}_{m_j}^{(j)}
\Bigl(x_j,x_{>j}-\sum_{k=j+1}^{n-1}m_k;a_j,a_{>j}+2\sum_{k=j+1}^{n-1}m_k\Bigr),
\label{Ridef}\\
\widetilde{\mathcal H}_i\mathcal{R}^{(i)}(\bm{x};\bm{a})&
=\sum_{k=i}^{n-1}m_k\Bigl(\sum_{k=i}^{n-1}m_k+\sum_{k=i}^na_k-1\Bigr)\mathcal{R}^{(i)}(\bm{x},\bm{a}).
\label{Rieig}\\
\widetilde{\mathcal H}_j\mathcal{R}^{(i)}(\bm{x};\bm{a})&
=0 \  \Leftrightarrow \ \widetilde{\mathcal H}_j\mathcal{R}^{(i)}(\bm{x};\bm{a})
=\mathcal{R}^{(i)}(\bm{x};\bm{a})\widetilde{\mathcal H}_j,\ \qquad j>i.
\end{align} 
The generalised forward and backward recursions hold for $\mathcal{R}^{(i)}$ as stated in the following
\begin{theo}
\label{theo:genfb}
\begin{align}
\text{Forward recursion:}\ \sum_{k=i}^{n}(x_k+a_k)\mathcal{R}^{(i)}(\bm{x}+\bm{e}_k)
&=\Bigl(\sum_{k=i}^n(x_k+a_k)+\sum_{k=i}^{n-1}m_k\Bigr)\mathcal{R}^{(i)}(\bm{x}),
\label{gfor}\\
\text{Backward recursion:}\qquad \ \ \sum_{k=i}^{n}x_k\mathcal{R}^{(i)}(\bm{x}-\bm{e}_k)
&=\Bigl(\sum_{k=i}^nx_k-\sum_{k=i}^{n-1}m_k\Bigr)\mathcal{R}^{(i)}(\bm{x}).
\label{gback}
\end{align}
\end{theo}
These can be shown by induction. At  the top, $i=n-1$, $\mathcal{R}^{(i)}$ is $\mathcal{P}^{(n-1)}$ and  it satisfies \eqref{gfor}
and \eqref{gback}  by \eqref{rec1} and  \eqref{rec2}. 
Suppose \eqref{gfor} and \eqref{gback} hold at $i$, the next level reads
\begin{align}
\mathcal{R}^{(i-1)}(\bm{x})=\mathcal{P}_{m_{i-1}}^{(i-1)}\bigl(x_{i-1},x_{>i-1}-\sum_{k=i}^{n-1}m_k;a_{i-1},a_{>i-1}+2\sum_{k=i}^{n-1}m_k\bigr)\mathcal{R}^{(i)}(\bm{x}),
\end{align}
in which $\mathcal{R}^{(i)}(\bm{x})$ does not contain $x_{i-1}$. Backward recursion comes first.
\begin{align*}
&\sum_{k=i-1}^nx_k\mathcal{R}^{(i-1)}(\bm{x}-\bm{e}_k)\\
&\qquad =x_{i-1}\mathcal{R}^{(i-1)}(\bm{x}-\bm{e}_{i-1})
+\sum_{k=i}^nx_k\mathcal{R}^{(i-1)}(\bm{x}-\bm{e}_k)\\
& \qquad =x_{i-1}\mathcal{P}_{m_{i-1}}^{(i-1)}\bigl(x_{i-1}-1,x_{>i-1}-\sum_{k=i}^{n-1}m_k;a_{i-1},a_{>i-1}+2\sum_{k=i}^{n-1}m_k\bigr)\mathcal{R}^{(i)}(\bm{x})\\
& \qquad \ +
\bigl(\sum_{k=i}^nx_k-\sum_{k=i}^{n-1}m_k\bigr)\mathcal{P}_{m_{i-1}}^{(i-1)}\bigl(x_{i-1},x_{>i-1}-1-\sum_{k=i}^{n-1}m_k;a_{i-1},a_{>i-1}+2\sum_{k=i}^{n-1}m_k\bigr)\mathcal{R}^{(i)}(\bm{x})\\
&\qquad =\bigl(\sum_{k=i-1}^nx_k-\sum_{k=i}^{n-1}m_k\bigr)\mathcal{P}_{m_{i-1}}^{(i-1)}\bigl(x_{i-1},x_{>i-1}-\sum_{k=i}^{n-1}m_k;a_{i-1},a_{>i-1}+2\sum_{k=i}^{n-1}m_k\bigr)\mathcal{R}^{(i)}(\bm{x})\\
&\qquad =\bigl(\sum_{k=i-1}^nx_k-\sum_{k=i-1}^{n-1}m_k\bigr)\mathcal{R}^{(i-1)}(\bm{x}).
\end{align*}
The assumption is used for the second term and \eqref{rec2} used for $\mathcal{P}_{m_{i-1}}^{(i-1)}$.
For  forward recursion \eqref{gfor} it goes similarly.

\begin{align*}
&\sum_{k=i-1}^n(x_k+a_k)\mathcal{R}^{(i-1)}(\bm{x}+\bm{e}_k)\\
&\qquad =(x_{i-1}+a_{i-1})\mathcal{R}^{(i-1)}(\bm{x}+\bm{e}_{i-1})
+\sum_{k=i}^n(x_k+a_k)\mathcal{R}^{(i-1)}(\bm{x}+\bm{e}_k)\\
& \qquad =(x_{i-1}+a_{i-1})\mathcal{P}_{m_{i-1}}^{(i-1)}\bigl(x_{i-1}+1,x_{>i-1}-\sum_{k=i}^{n-1}m_k;a_{i-1},a_{>i-1}+2\sum_{k=i}^{n-1}m_k\bigr)\mathcal{R}^{(i)}(\bm{x})\\
& \qquad \ +
\bigl(\sum_{k=i}^n(x_k+a_k)\!+\!\!\sum_{k=i}^{n-1}m_k\bigr)
\mathcal{P}_{m_{i-1}}^{(i-1)}
\bigl(x_{i-1},x_{>i-1}+\!1-\!\sum_{k=i}^{n-1}m_k;a_{i-1},a_{>i-1}\!+\!2\sum_{k=i}^{n-1}m_k\bigr)
\mathcal{R}^{(i)}(\bm{x})\\
&\qquad =\bigl(\sum_{k=i-1}^n(x_k+a_k)+\!\!\sum_{k=i-1}^{n-1}m_k\bigr)
\mathcal{P}_{m_{i-1}}^{(i-1)}\bigl(x_{i-1},x_{>i-1}-\!\sum_{k=i}^{n-1}m_k;a_{i-1},a_{>i-1}+2\sum_{k=i}^{n-1}m_k\bigr)\mathcal{R}^{(i)}(\bm{x})\\
&\qquad =\bigl(\sum_{k=i-1}^n(x_k+a_k)+\sum_{k=i-1}^{n-1}m_k\bigr)\mathcal{R}^{(i-1)}(\bm{x}).
\end{align*}
The assumption is used for the second term and \eqref{rec1} 
used for $\mathcal{P}_{m_{i-1}}^{(i-1)}$. \hfill$\square$\\
By using these results, the {\bf Eigenvalue Theorem \ref{theo:eigens}} for $\mathcal{P}^{(i)}$ is rewritten
for $\mathcal{R}^{(i)}$. The proof in Appendix B goes almost parallel for $\mathcal{R}^{(i)}$. 
\begin{theo}
\label{theo:eigensR}
{\bf Eigenvalue Theorem}
\begin{align}
\widetilde{\mathcal H}_T\mathcal{R}^{(i)}\bigl(\bm{x};\bm{a}\bigr)&
=\sum_{k=i}^{n-1}m_k\Bigl(\sum_{k=i}^{n-1}m_k+|a|+b-1\Bigr)\mathcal{R}^{(i)}\bigl(\bm{x};\bm{a}\bigr),
\label{eigTG}\\
\widetilde{\mathcal H}_j\mathcal{R}^{(i)}\bigl(\bm{x};\bm{a}\bigr)&
=\sum_{k=j}^{n-1}m_k\Bigl(\sum_{k=j}^{n-1}m_k+a_j+\cdots+a_n-1\Bigr)
\mathcal{R}^{(i)}\bigl(\bm{x};\bm{a}\bigr),\quad j<i,
\label{eigjG}\\
&=\sum_{k=j}^{n-1}m_k\Bigl(\sum_{k=j}^{n-1}m_k+a_j+a_{>j}-1\Bigr)
\mathcal{R}^{(i)}\bigl(\bm{x};\bm{a}\bigr),\qquad \quad j<i,\n
\widetilde{\mathcal H}_i\mathcal{R}^{(i)}\bigl(\bm{x};\bm{a}\bigr)&
=\sum_{k=i}^{n-1}m_k\Bigl(\sum_{k=i}^{n-1}m_k+a_i+\cdots+a_n-1\Bigr)
\mathcal{R}^{(i)}\bigl(\bm{x};\bm{a}\bigr),
\label{eigiG}\\
&=\sum_{k=i}^{n-1}m_k\Bigl(\sum_{k=i}^{n-1}m_k+a_i+a_{>i}-1\Bigr)
\mathcal{R}^{(i)}\bigl(\bm{x};\bm{a}\bigr),\n
\widetilde{\mathcal H}_j\mathcal{R}^{(i)}\bigl(\bm{x};\bm{a}\bigr)&=0,\quad j>i,
\hspace{50mm}
i=1,\ldots,n-1.
\label{eig0G}
\end{align}
\end{theo}
The final step is to combine $R^{(1)}$ with a polynomial $\psi$ in $|x|$, in such a way that 
$R^{(1)}\psi$ is the eigenpolynomial of $\widetilde{\mathcal H}_T$.
As $\widetilde{\mathcal H}_T=\widetilde{\mathcal H}_0+\widetilde{\mathcal H}_1$ and
$\widetilde{\mathcal H}_1=\sum_{j\neq k}x_j(x_k+a_k)(1-e^{-\partial_j+\partial_k})$ 
and $e^{-\partial_j+\partial_k}\psi(|x|)=\psi(|x|)$,
\begin{align*}
\widetilde{\mathcal H}_1\mathcal{R}^{(1)}\psi(|x|)
=\Bigl(\widetilde{\mathcal H}_1\mathcal{R}^{(1)}\Bigr)\psi(|x|)=
\Bigl(\sum_{j=1}^{n-1}m_j\Bigr)\Bigl(\sum_{j=1}^{n-1}m_j+|a|-1\Bigr)\mathcal{R}^{(1)}\psi(|x|).
\end{align*}
For $\widetilde{\mathcal H}_0
=\sum_{j=1}^n(N-|x|)(x_k+a_k)(1-e^{\partial_j})+\sum_{j=1}^nx_j(N-|x|+b)(1-e^{-\partial_j})$,
\begin{align*}
&\widetilde{\mathcal H}_0\mathcal{R}^{(1)}\psi(|x|)\n
&\ =(N-|x|)(|x|+|a|)\mathcal{R}^{(1)}\psi(|x|)
 -(N-|x|)\sum_{j=1}^n(x_j+a_j)\mathcal{R}^{(1)}(\bm{x}+\bm{e}_j)\psi(|x|+1)\n
&\quad +|x|(N-|x|+b)\mathcal{R}^{(1)}\psi(|x|)-(N-|x|+b)
\sum_{j=1}^nx_j\mathcal{R}^{(1)}(\bm{x}-\bm{e}_j)\psi(|x|-1)\n
&\ =\mathcal{R}^{(1)}\left\{\Bigl((N-|x|)(|x|+|a|)+|x|(N-|x|+b)\Bigr)\psi(|x|)\right.\n
&\qquad \qquad -\big(N-|x|\big)\Bigl(|x|+|a|+\sum_{j=1}^{n-1}m_j\Bigr)\psi(|x|+1)\n
&\qquad \qquad -\left.\Bigl(|x|-\sum_{j=1}^{n-1}m_j\Bigr)\bigl(N-|x|+b\bigr)\psi(|x|-1)\right\}.
\end{align*}
Here {\bf Theorem \ref{theo:genfb}} is used. Finally this is rearranged to
\begin{align*}
&\widetilde{\mathcal H}_0\mathcal{R}^{(1)}\psi(|x|)\n
&\ =\mathcal{R}^{(1)}\left\{\Bigl(N-\sum_{j=1}^{n-1}m_j-\bigl(|x|-\sum_{j=1}^{n-1}m_j\bigr)\Bigr)
\Bigl(|x|-\sum_{j=1}^{n-1}m_j+|a|+2\sum_{j=1}^{n-1}m_j\Bigr)\right.\n
&\hspace{40mm}\times\Bigl(\psi(|x|)-\psi(|x|+1)\Bigr)\n
& \hspace{15mm} +\Bigl(|x|-\sum_{j=1}^{n-1}m_j\Bigr)
\Big(N-\sum_{j=1}^{n-1}m_j-\bigl(|x|-\sum_{j=1}^{n-1}m_j\bigr)+b\Bigr)\n
&\hspace{40mm} \left. \times\Bigl(\psi(|x|)-\psi(|x|-1)\Bigr)\right\}
\quad  +b\Bigl(\sum_{j=1}^{n-1}m_j\Bigr)\mathcal{R}^{(1)}\psi(|x|).
\end{align*}
The main part is the difference equation for the single variable Hahn polynomial \eqref{Hahneq} with the shifts
\begin{equation*}
x \to |x|-\sum_{j=1}^{n-1}m_j,\quad a \to |a|+2\sum_{j=1}^{n-1}m_j,\quad N \to N-\sum_{j=1}^{n-1}m_j,
\end{equation*}
with a degree $m_0$ leading to 
\begin{align*}
{}_3F_2\left(
  \genfrac{}{}{0pt}{}{-m_0,\,m_0+2\sum_{j=1}^{n-1}m_j+|a|+b-1,\,-|x|+\sum_{j=1}^{n-1}m_j}
  {|a|+2\sum_{j=1}^{n-1}m_j,\,-N+\sum_{j=1}^{n-1}m_j}\Bigm|1\right),
\tag{\ref{sol2}}
\end{align*}
with the eigenvalue
\begin{align*}
\mathcal{E}_T(\bm{m})&=\sum_{j=1}^{n-1}m_j\Bigl(\sum_{j=1}^{n-1}m_j+|a|-1\Bigr)+m_0\Bigl(m_0+|a|+b-1+2\sum_{j=1}^{n-1}m_j\Bigr)+
b\sum_{j=1}^{n-1}m_j\\
&=|m|\Bigl(|m|+|a|+b-1\Bigr),\quad  \bm{m}=(m_0,m_1,\ldots,m_{n-1}) \tag{\ref{Htieq}}.
\end{align*}

\section{Multivariate Krawtchouk polynomials}
\label{sec:newmKra}
It is known that the single variable Krawtchouk polynomials are obtained from the single variable
Hahn polynomials \eqref{Hahnpoly} by sending $a$ and $b$ to infinity by the same rate \cite{KLS}(9.5.16).
Let us define multivariate Krawtchouk polynomials from the multivariate Hahn polynomials \eqref{sol1}
by the limits $a_i\to a_ib$ and $b\to \infty$. 
The polynomials are defined on the same  compact $n$-dimensional
integer lattice $\cX$ \eqref{XKdef} as that of the Hahn polynomials.
Most Propositions and Theorems in section \ref{sec:Problem} to \ref{sec:geneig} hold in the  limit 
and simplified except for
the definition of $\mathcal{P}_m^{(i)}$ \eqref{Piform} which needs to be  divided by $(a_i)_m$  to make it finite.
Since the parameters $\{a_i\}$ are scaled, they do not change in the shape invariance calculations, 
like  the parameter $p$ in the single variable Krawtchouk case \cite{KLS,os12}.
Here I present the basic results of the multivariate  Krawtchouk polynomials mostly without proof or derivation.

The difference operators governing the  multivariate  Krawtchouk polynomials  are
\label{def:Ht0K}
\begin{align}
\widetilde{\mathcal H}_{KT}&=\widetilde{\mathcal H}_{K0}+\widetilde{\mathcal H}_{K1},
\label{HtsumK}\\
\widetilde{\mathcal H}_{K0}&\eqdef \sum_{j=1}^n\Bigl((N-|x|)a_j(1-e^{\partial_j})+x_j(1-e^{-\partial_j})\Bigr)
=\sum_{j=1}^n\bigl((N-|x|)a_j-x_je^{-\partial_j}\bigr)\bigl(1-e^{\partial_j}\bigr),
\label{Ht0defK}\\
\widetilde{\mathcal H}_{K1}&\eqdef \sum_{j\neq k\ge1}^n x_ja_k\bigl(1-e^{-\partial_j+\partial_k}\bigr),
\quad \widetilde{\mathcal H}_{Ki}\eqdef \sum_{j\neq k\ge i}^n x_ja_k\bigl(1-e^{-\partial_j+\partial_k}\bigr),
\label{Ht1defK}\\
&[\widetilde{\mathcal H}_{KT},\widetilde{\mathcal H}_{Ki}]
=[\widetilde{\mathcal H}_{Ki},\widetilde{\mathcal H}_{Kj}]=0.\qquad i,j=0,1,\ldots,n.
\end{align}
The corresponding orthogonality measure is the  multinomial distribution depending on $\bm{a}$,
\begin{align}
&W_K(\bm{x};\bm{a},N)=\frac{N!}{x_1!\cdots x_n!x_0!}
\frac{\prod_{i=1}^na_i^{x_i}}{(1+|a|)^N},
\quad \bm{a}=(a_1,\ldots,a_n)\in\mathbb{R}_{>0}^n.
\label{WK}
\end{align} 

As in the multi-Hahn case, all the eigenvalues of $\widetilde{\mathcal H}_{KT}$ are the same, as 
\begin{align*}
\widetilde{\mathcal H}_{KT}p_m(x_J)&=(N-x_J)a_J\bigl(p_m(x_J)-p_m(x_J+1)\bigr)
+x_J(1+|a|-a_J)\bigl(p_m(x_J)-p_m(x_J-1)\bigr)\\
&=m(|a|+1)p_m(x_J),\  p_m(x)=K_m(x;p,N)={}_2F_1\Bigl(
  \genfrac{}{}{0pt}{}{-m,\,-x}{-N}\Bigm|p^{-1}\Bigr), p=\frac{a_J}{1+|a|},
\end{align*}
in which $x_J$ is defined in \eqref{xjdef} and $K_m$ is the single variable Krawtchouk polynomial 
\cite{KLS,os12} of degree $m$.
The type two eigenpolynomials are
\begin{align}
 \mathcal{P}_{Km}^{(i)}\bigl(\bm{x})&= \mathcal{P}_{Km}^{(i)}\bigl(x_i,x_{>i};a_i,a_{>i}\bigr)
=\sum_{k=0}^{m}(-1)^k\binom{m}{k}\left(\frac{a_{>i}}{a_i}\right)^k(-x_i)_{k}(-x_{>i})_{m-k},
 \label{PiformK}\\
 \widetilde{\mathcal H}_{KT}\mathcal{P}_{Km}^{(i)}\bigl(\bm{x})
&=m(|a|+1)\mathcal{P}_{Km}^{(i)}(\bm{x}),
\label{eigHTK}\\
\widetilde{\mathcal H}_{Ki }\mathcal{P}_{Km}^{(i)}\bigl(\bm{x})
&
=m(a_i+a_{>i})\mathcal{P}_{Km}^{(i)}(\bm{x}),
\label{eigHiK}\\
\widetilde{\mathcal H}_{Kj }\mathcal{P}_{Km}^{(i)}\bigl(\bm{x})
&=m(a_j+a_{>j})\mathcal{P}_{Km}^{(i)}(\bm{x}), \qquad j<i,\n
\widetilde{\mathcal H}_{Kj }\mathcal{P}_{Km}^{(i)}\bigl(\bm{x})&=0,\quad j>i
\quad \Rightarrow (\mathcal{P}_{Km}^{(i)},\mathcal{P}_{Km}^{(j)})=0,\quad i\neq j.
\label{PmorthoK}\\
&\mathcal{P}_{Km}^{(i)}\bigl(0,m;a_i,a_{>i})=(-1)^mm!.
\label{sPvalK}
\end{align}
$\mathcal{P}_{Km}^{(i)}\bigl(\bm{x})$ is obtained by dividing 
$\mathcal{P}_{m}^{(i)}\bigl(\bm{x};\bm{\lambda})$ \eqref{Piform} by $(a_i)_m$ and $a_j \to a_jt$, $t\to\infty$
and the summation index $k$ is changed to $m-k$, for simpler appearance.
 \begin{align}
\text{Forward shift:}\quad &\mathcal{P}_{Km}^{(i)}\bigl(x_i,x_{>i}+1;a_i,a_{>i}\bigr)
-\mathcal{P}_{Km}^{(i)}\bigl(x_i+1,x_{>i};a_i,a_{>i}\bigr)\n
&\quad  =-\frac{m(a_i+a_{>i})}{a_i}\mathcal{P}_{Km-1}^{(i)}\bigl(x_i,x_{>i};a_i,a_{>i}\bigr),
 \label{PiforK}\\
\text{Backward shift:}\quad   &x_{>i}a_i\mathcal{P}_{Km}^{(i)}\bigl(x_i,x_{>i}-1;a_i,a_{>i}\bigr)
 -x_ia_{>i}\mathcal{P}_{Km}^{(i)}\bigl(x_i-1,x_{>i};a_i,a_{>i}\bigr)\n
&\quad =-a_i\mathcal{P}_{Km+1}^{(i)}\bigl(x_i,x_{>i};a_i,a_{>i}\bigr),
\label{PibackK}
 \end{align}
 in which the sign change  and the factor $a_i$ are  due to the redefinition of 
 $\mathcal{P}_{Km}^{(i)}$ \eqref{PiformK}. 
  \begin{align}
 \text{Forward recursion:}\quad &a_{i}\mathcal{P}_{Km}^{(i)}\bigl(x_i+1,x_{>i};a_i,a_{>i}\bigr)
                                                           +a_{>i}\mathcal{P}_{Km}^{(i)}\bigl(x_i,x_{>i}+1;a_i,a_{>i}\bigr)\n                                                        
 &\quad =(a_i+a_{>i}) \mathcal{P}_{Km}^{(i)}\bigl(x_i,x_{>i};a_i,a_{>i}\bigr),
\label{rec1K}\\
 \text{Backward recursion:}\quad & x_{i}\mathcal{P}_{Km}^{(i)}\bigl(x_i-1,x_{>i};a_i,a_{>i}\bigr)
                                            +x_{>i}\mathcal{P}_{Km}^{(i)}\bigl(x_i,x_{>i}-1;a_i,a_{>i}\bigr)\n  
&\quad =(x_i+x_{>i}-m) \mathcal{P}_{Km}^{(i)}\bigl(x_i,x_{>i};a_i,a_{>i}\bigr).
\label{rec2K}
\end{align}  
Starting from $\mathcal{P}_K^{(n-1)}$, the combination goes as far as  $\mathcal{P}_K^{(i)}$ to produce
\begin{align}
\mathcal{R}_K^{(i)}(\bm{x};\bm{a})&=\prod_{j=i}^{n-1}\mathcal{P}_{Km_j}^{(j)}
\Bigl(x_j,x_{>j}-\sum_{k=j+1}^{n-1}m_k;a_j,a_{>j}\Bigr),
\label{RidefK}\\
\widetilde{\mathcal H}_{Ki}\mathcal{R}_K^{(i)}(\bm{x};\bm{a})&
=\sum_{k=i}^{n-1}m_k\Bigl(\sum_{k=i}^na_k\Bigr)\mathcal{R}_K^{(i)}(\bm{x};\bm{a}),
\label{RieigK}\\
\widetilde{\mathcal H}_{Kj}\mathcal{R}_K^{(i)}(\bm{x};\bm{a})&
=0 \  \Leftrightarrow \ \widetilde{\mathcal H}_{Kj}\mathcal{R}_K^{(i)}(\bm{x};\bm{a})
=\mathcal{R}_K^{(i)}(\bm{x};\bm{a})\widetilde{\mathcal H}_{Kj},\ \qquad j>i.
\end{align}  
\begin{align}
\text{Forward recursion:}\quad \ \sum_{k=i}^{n}a_k\mathcal{R}_K^{(i)}(\bm{x}+\bm{e}_k;\bm{a})
&=\Bigl(\sum_{k=i}^na_k\Bigr)\mathcal{R}_K^{(i)}(\bm{x};\bm{a}),
\label{gforK}\\
\text{Backward recursion:}\quad \ \ \sum_{k=i}^{n}x_k\mathcal{R}_K^{(i)}(\bm{x}-\bm{e}_k;\bm{a})
&=\Bigl(\sum_{k=i}^nx_k-\sum_{k=i}^{n-1}m_k\Bigr)\mathcal{R}_K^{(i)}(\bm{x};\bm{a}).
\label{gbackK}
\end{align}
\begin{theo}
\label{theo:mainK}{\bf Main Theorem (K)}
The complete set of eigenpolynomials of the multivariate Krawtchouk system is
\begin{align}
\mathcal{P}_{K\bm m}(\bm{x};\bm{a},N)
&=\prod_{i=1}^{n-1}\mathcal{P}_{Km_i}^{(i)}
\Bigl(x_i,x_{>i}\!-\!\sum_{j=i+1}^{n-1}m_j;a_i,a_{>i}\Bigr)
 K_{m_0}\Bigl(|x|-\!\sum_{j=1}^{n-1}m_j;\frac{|a|}{|a|+1},N\!-\!\sum_{j=1}^{n-1}m_j\Bigr)
\label{sol1K}\\
&=\prod_{i=1}^{n-1} \sum_{k=0}^{m_i}(-1)^k\binom{m_i}{k}\Bigl(\frac{a_{>i}}{a_i}\Bigr)^k
 \times\bigl(-x_i\bigr)_{k}\Bigl(-x_{>i}+\sum_{j=i+1}^{n-1}m_j\Bigr)_{m-k}\n
&\qquad \times {}_2F_1\left(
  \genfrac{}{}{0pt}{}{-m_0,\,-|x|+\sum_{j=1}^{n-1}m_j}
  {-N+\sum_{j=1}^{n-1}m_j}\Bigm|\frac{|a|+1}{|a|}\right),
  \label{sol2K}\\[4pt]
& \hspace{-20mm}\widetilde{\mathcal H}_{KT}\mathcal{P}_{K\bm m}(\bm{x};\bm{a},N)
=|m|(|a|+1)\mathcal{P}_{K\bm m}(\bm{x};\bm{a},N),
\label{HteqK}\\
& \hspace{-20mm}\widetilde{\mathcal H}_{Ki}\mathcal{P}_{K\bm m}(\bm{x};\bm{a},N)
=|m|\Bigl(\sum_{j=i}^{n}a_j\Bigr)\mathcal{P}_{K\bm m}(\bm{x};\bm{a},N),\qquad i=1,\ldots, n-1,
\label{HtieqK}\\
&\hspace{20mm}  \bm{m}=(m_0,m_1,\ldots,m_{n-1})\in\cX,\n
&\Big(\mathcal{P}_{K\bm m},\mathcal{P}_{K\bm{m}'}\Bigr)=0,\qquad \bm{m}\neq\bm{m}'.
\label{genorthoK}
\end{align}
\end{theo}  

\section{Multivariate Meixner polynomials}
\label{sec:newmMei}
Let us define multivariate Meixner polynomials from the multivariate Hahn polynomials 
by the parameter redefinition and another limiting procedure \cite{tra3}, {\em i.e.\/},
\begin{equation}
N\to -\beta,\quad \beta>0,\quad a_j\to -a_jb,\quad a_j>0,\quad j=1,\ldots,n,\quad |a|<1,\quad b\to+\infty.
\label{Meilimit}
\end{equation}
The polynomials are now defined on a semi-infinite integer lattice $\cX=\mathbb{N}_0^n$. 
The difference operators governing the  multivariate  Meixner polynomials  are
\label{def:Ht0M}
\begin{align}
\widetilde{\mathcal H}_{MT}&=\widetilde{\mathcal H}_{M0}+\widetilde{\mathcal H}_{M1},
\label{HtsumM}\\
\widetilde{\mathcal H}_{M0}&\eqdef \sum_{j=1}^n\Bigl((\beta+|x|)a_j(1-e^{\partial_j})+x_j(1-e^{-\partial_j})\Bigr)
=\sum_{j=1}^n\bigl((\beta+|x|)a_j-x_je^{-\partial_j}\bigr)\bigl(1-e^{\partial_j}\bigr),
\label{Ht0defM}\\
\widetilde{\mathcal H}_{M1}&\eqdef -\sum_{j\neq k\ge1}^n x_ja_k\bigl(1-e^{-\partial_j+\partial_k}\bigr)
=-\widetilde{\mathcal H}_{K1},
\quad \widetilde{\mathcal H}_{Mi}\eqdef -\widetilde{\mathcal H}_{Ki},
\label{Ht1defM}\\
&[\widetilde{\mathcal H}_{MT},\widetilde{\mathcal H}_{Mi}]
=[\widetilde{\mathcal H}_{Mi},\widetilde{\mathcal H}_{Mj}]=0.\qquad i,j=0,1,\ldots,n.
\end{align}
The corresponding orthogonality measure is  the negative multinomial distribution {\rm \cite{Gri12}}, 
\begin{equation}
W_M({\bm x};\bm{a},\beta)\eqdef \frac{(\beta)_{|x|}\bm{a}^{\bm x}}{\bm{x}!}\bigl(1-|a|\bigr)^{\beta},
\quad |a|<1.
\label{nMeiWdef}
\end{equation}

As in the multi-Hahn case, all the eigenvalues of $\widetilde{\mathcal H}_{MT}$ are the same, as 
\begin{align*}
\widetilde{\mathcal H}_{MT}p_m(x_J)&=(\beta+x_J)a_J\bigl(p_m(x_J)-p_m(x_J+1)\bigr)
+x_J(1-|a|+a_J)\bigl(p_m(x_J)-p_m(x_J-1)\bigr)\\
&=m\bigl(1-|a|\bigr)p_m(x_J),\\
  p_m(x)&=M_m(x;c,\beta)={}_2F_1\Bigl(
  \genfrac{}{}{0pt}{}{-m,\,-x}{\beta}\Bigm|1-c^{-1}\Bigr), \quad c=\frac{a_J}{1-|a|+a_J},
\end{align*}
in which $x_J$ is defined in \eqref{xjdef} and $M_m$ is the single variable 
Meixner polynomial \cite{KLS,os12} of degree $m$.
Since $\widetilde{\mathcal H}_{Mi}=-\widetilde{\mathcal H}_{Ki}$, $i=1,\ldots,n$, all the type two
eigenpolynomials are the same as those of the multivariate Krawtchouk polynomials with negative eigenvalues,
\begin{align}
 \mathcal{P}_{Mm}^{(i)}\bigl(\bm{x})&= \mathcal{P}_{Mm}^{(i)}\bigl(x_i,x_{>i};a_i,a_{>i}\bigr)
=\sum_{k=0}^{m}(-1)^k\binom{m}{k}\left(\frac{a_{>i}}{a_i}\right)^k(-x_i)_{k}(-x_{>i})_{m-k},
 \label{PiformM}\\
 \widetilde{\mathcal H}_{MT}\mathcal{P}_{Mm}^{(i)}\bigl(\bm{x})
&=m(1-|a|)\mathcal{P}_{Mm}^{(i)}(\bm{x}),
\label{eigHTM}\\
\widetilde{\mathcal H}_{Mi }\mathcal{P}_{Mm}^{(i)}\bigl(\bm{x})
&
=-m(a_i+a_{>i})\mathcal{P}_{Mm}^{(i)}(\bm{x}),
\label{eigHiM}\\
\widetilde{\mathcal H}_{Mj }\mathcal{P}_{Mm}^{(i)}\bigl(\bm{x})
&=-m(a_j+a_{>j})\mathcal{P}_{Mm}^{(i)}(\bm{x}), \qquad j<i,\n
\widetilde{\mathcal H}_{Mj }\mathcal{P}_{Mm}^{(i)}\bigl(\bm{x})&=0,\quad j>i
\quad \Rightarrow (\mathcal{P}_{Mm}^{(i)},\mathcal{P}_{Mm}^{(j)})=0,\quad i\neq j.
\label{PmorthoM}\\
&\mathcal{P}_{Mm}^{(i)}\bigl(0,m;a_i,a_{>i})=(-1)^mm!.
\label{sPvalM}
\end{align}
The propositions and theorems related to the type two eigenpolynomials are the same as those given for the
multi-Krawtchouk polynomials, except for the opposite signs of the eigenvalues.
Therefore I present only the following 
\begin{theo}
\label{theo:mainM}{\bf Main Theorem (M)}
The complete set of eigenpolynomials of the multivariate Meixner system is
\begin{align}
\mathcal{P}_{M\bm m}(\bm{x};\bm{a},\beta)
&=\prod_{i=1}^{n-1}\mathcal{P}_{Mm_i}^{(i)}
\Bigl(x_i,x_{>i}\!-\!\sum_{j=i+1}^{n-1}m_j;a_i,a_{>i}\Bigr)
 M_{m_0}\Bigl(|x|-\!\sum_{j=1}^{n-1}m_j;|a|,\beta\!+\!\sum_{j=1}^{n-1}m_j\Bigr)
\label{sol1M}\\
&=\prod_{i=1}^{n-1} \sum_{k=0}^{m_i}(-1)^k\binom{m_i}{k}\Bigl(\frac{a_{>i}}{a_i}\Bigr)^k
 \times\bigl(-x_i\bigr)_{k}\Bigl(-x_{>i}+\sum_{j=i+1}^{n-1}m_j\Bigr)_{m-k}\n
&\qquad \times {}_2F_1\left(
  \genfrac{}{}{0pt}{}{-m_0,\,-|x|+\sum_{j=1}^{n-1}m_j}
  {\beta+\sum_{j=1}^{n-1}m_j}\Bigm|\frac{|a|-1}{|a|}\right),
  \label{sol2M}\\[4pt]
& \hspace{-20mm}\widetilde{\mathcal H}_{MT}\mathcal{P}_{M\bm m}(\bm{x};\bm{a},\beta)
=|m|\bigl(1-|a|\bigr)\mathcal{P}_{M\bm m}(\bm{x};\bm{a},\beta),
\label{HteqM}\\
& \hspace{-20mm}\widetilde{\mathcal H}_{Mi}\mathcal{P}_{M\bm m}(\bm{x};\bm{a},\beta)
=-|m|\Bigl(\sum_{j=i}^{n}a_j\Bigr)\mathcal{P}_{M\bm m}(\bm{x};\bm{a},\beta),\qquad i=1,\ldots, n-1,
\label{HtieqM}\\
&\hspace{20mm}  \bm{m}=(m_0,m_1,\ldots,m_{n-1})\in\cX,\n
&\Big(\mathcal{P}_{M\bm m},\mathcal{P}_{M\bm{m}'}\Bigr)=0,\qquad \bm{m}\neq\bm{m}'.
\label{genorthoM}
\end{align}
\end{theo}  
\section{Comments}
\label{sec:comments}
Some comments are in order.\\
$\bullet$ Hypergeometric classification. It would be quite interesting to see how the obtained
eigenpolynomials 
$\mathcal{P}_{\bm m}(\bm{x};\bm{a},b,N)$ \eqref{sol1}, $\mathcal{P}_{K\bm m}(\bm{x};\bm{a},N)$ \eqref{sol1K}
and $\mathcal{P}_{M\bm m}(\bm{x};\bm{a},\beta)$ \eqref{sol1M}
 fit in the general scheme of hypergeometric functions
of Aomoto-Gelfabd \cite{AK,gelfand}.\\
$\bullet$ $\mathfrak{S}_n$ symmetry. The difference operator $\widetilde{\mathcal H}_T$ \eqref{HTdef}
is $\mathfrak{S}_n$ symmetric $\bm{x}\to\sigma\bm{x}$, $\bm{a}\to\sigma\bm{a}$,  
$\forall\sigma\in\mathfrak{S}_n$. 
But the eigenpolynomials 
$\mathcal{P}_{\bm m}(\bm{x};\bm{a},b,N)$ \eqref{sol1} are not  $\mathfrak{S}_n$ symmetric, as 
a very special order \eqref{tidef} is chosen for orthogonalising  eigenpolynomials.
Applying an arbitrary $\mathfrak{S}_n$ transformation would make the eigenpolynomials quite clumsy.
 A passable change would be the mirror
reflection
\begin{equation}
(x_1,\ldots,x_{n-1},x_n) \to (x_n,x_{n-1},\ldots,x_1),\qquad 
(a_1,\ldots,a_{n-1},a_n) \to (a_n,a_{n-1},\ldots,a_1).
\end{equation}
$\bullet$ Normalisation and norm. I do not know the principle to fix the normalisation of the 
eigenpolynomials  $\mathcal{P}_{\bm m}(\bm{x};\bm{a},b,N)$ \eqref{sol1}, 
$\mathcal{P}_{K\bm m}(\bm{x};\bm{a},N)$ \eqref{sol1K}
and $\mathcal{P}_{M\bm m}(\bm{x};\bm{a},\beta)$ \eqref{sol1M}. 
Thus their norms are not calculated.\\
$\bullet$ Birth and Death process. The stochastic process that Karlin-McGrgor \cite{KarMcG} intended to solve
is now solved exactly and explicitly, once the norms of the eigenpolynomials 
$\mathcal{P}_{\bm m}(\bm{x};\bm{a},b,N)$ \eqref{sol1} are calculated, as shown in \cite{bdsol}, \cite{mKrawt}\S4.7.\\
$\bullet$ Relation with many seminal works.\\
I have not yet deciphered the work of Karlin-McGregor \cite{KarMcG},  as it aims the solution 
of the time-dependent birth-death like problem, not the eigenvalue problem.\\
The eigenpolynomials 
$\mathcal{P}_{\bm m}(\bm{x};\bm{a},b,N)$ \eqref{sol1}, 
 look similar to those given by Tratnik \cite{tra3},
which are derived by a certain limit from the proposed multivariate Racah polynomials. 
The corresponding Krawtchouk polynomials $\mathcal{P}_{K\bm m}(\bm{x};\bm{a},N)$ \eqref{sol1K} 
and Meixner polynomials $\mathcal{P}_{M\bm m}(\bm{x};\bm{a},\beta)$ \eqref{sol1M} 
also look quite similar to those in Tratnik's \cite{tra3}.
Difference equation is not imposed in Tratnik's work \cite{tra3}.\\
Similar expressions are also found in the works of 
Iliev-Xu \cite{ilxu}(5.30), Xu \cite{xu3}(2.18),  as eigenvectors  of a difference operator, 
which is slightly different from mine 
$\widetilde{\mathcal H}_T$ \eqref{HTdef}.\\
 For $n=2$, my results $\mathcal{P}_{\bm m}(\bm{x};\bm{a},b,N)$ \eqref{sol1} are the same as those given by Genest-Vinet \cite{genest}A.1., including two difference operators
 and the forward-backward relations \cite{genest}A.6.
 Probably, their two variable generating function could be generalised to $n$-variable one. 
Since their method of derivation, the  overlapping coefficients \cite{genest},  
provides exactly the same results with the difference equation method for $n=2$, it seems possible to make bridges
connecting various methods, intertwining functions \cite{dunkl0,scarabotti}, 
coupling coefficients \cite{rosengren}.\\
The solution method of the present work could be deemed as  separation of the variables as
Karlin-McGregor say\cite{KarMcG}. In the eigenpolynomials 
$\mathcal{P}_{\bm m}(\bm{x};\bm{a},b,N)$ \eqref{sol1}, the last part $H_{m_0}(|x|)$ corresponds to the radial part and
the rest are the spherical harmonics.
\section*{Acknowledgements}
RS thanks Mourad Ismail for inducing him to explore the mysterious maze of birth and death processes.
He also thanks Kazuhiko Aomoto and Charles Dunkl for useful comments and warm encouragements.

\section*{Declarations}
\begin{itemize}
\item Funding: No funds, grants, or other support was received.
\item Data availability statement: Data sharing not applicable to this article as 
no datasets were generated or analysed during the current study.
\item Competing Interests: The author has no competing interests to declare that 
are relevant to the content of this article.
\end{itemize}

\appendix
\section*{Appendix A: Inductive proof of type two eigenpolynomials}
\label{appendA}
\setcounter{equation}{0}
\renewcommand{\theequation}{A.\arabic{equation}}

Assuming the formula \eqref{Piform} is correct at degree $m-1$, 
\begin{equation*}
\mathcal{P}_{m-1}^{(i)}(x_i,x_{>i};a_i+1,a_{>i}+1)\phi_0^{(i)}(\bm{x};\bm{\lambda}+\bm{\delta})
\propto \mathcal{B}(\bm{\lambda}+\bm{\delta})^T\cdots
  \mathcal{B}(\bm{\lambda}+(m-1)\bm{\delta})^T
  \phi_0(x\,;\bm{\lambda}+m\bm{\delta}),
\end{equation*}
I derive degree $m$ formula by applying 
$\mathcal{B}(\bm{\lambda})^T= \sqrt{x_{>i}(x_i+a_i)}\,e^{-\partial_{>i}}-\sqrt{x_i(x_{>i}+a_{>i})}\,e^{-\partial_i}$,
\begin{equation*}
\mathcal{P}_{m}^{(i)}(x_i,x_{>i};a_i,a_{>i})\phi_0^{(i)}(\bm{x};\bm{\lambda})
\propto\mathcal{B}(\bm{\lambda})^T\mathcal{P}_{m-1}^{(i)}(x_i,x_{>i};a_i+1,a_{>i}+1)
\phi_0^{(i)}(\bm{x};\bm{\lambda}+\bm{\delta}).
\end{equation*}
The part $(b')_{N-1-x_i-x_{>i}}/(N-1-x_i-x_{>i})!$ in $\phi_0^{(i)}(\bm{x};\bm{\lambda}+\bm{\delta})$ goes
to $(b')_{N-x_i-x_{>i}}/(N-x_i-x_{>i})!$ in $\phi_0^{(i)}(\bm{x};\bm{\lambda})$ by 
$e^{-\partial_{>i}}$ and $e^{-\partial_{i}}$ in $\mathcal{B}(\bm{\lambda})^T$. The constant parts, $(N-1)!/(|a|+b)_{N-1}$
in $\phi_0^{(i)}(\bm{x};\bm{\lambda}+\bm{\delta})$ and $N!/(|a|+b)_{N}$ in $\phi_0^{(i)}(\bm{x};\bm{\lambda})$
are neglected since they are immaterial. 
The essential part of the degree $m-1$ formula is
\begin{align}
&\mathcal{P}_{m-1}^{(i)}\bigl(x_i,x_{>i};a_i+1,a_{>i}+1\bigr)\phi_0^{(i)}(\bm{x};\bm{\lambda}+\bm{\delta})\n
&\qquad=
 \sum_{k=0}^{m-1}(-1)^k\binom{m-1}{k}(a_{>i}+1+k)_{m-1-k}(a_i+m-k)_k\times X,\n
 &\qquad \qquad  X\eqdef (-x_i)_{m-1-k}(-x_{>i})_k \sqrt{\frac{(a_i+1)_{x_i}(a_{>i}+1)_{x_{>i}}}{x_i!x_{>i}!}}.
\end{align}
I apply $\mathcal{B}(\bm{\lambda})^T$ to the variable dependent part $X$.
\begin{align}
&\mathcal{B}(\bm{\lambda})^TX\n
&\quad=(-x_i)_{m-1-k}(-x_{>i}+1)_k \sqrt{\frac{x_{>i}(x_i+a_i)(a_i+1)_{x_i}(a_{>i}+1)_{x_{>i}-1}}{x_i!(x_{>i}-1)!}}\n
&\quad\ -(-x_i+1)_{m-1-k}(-x_{>i})_k
\sqrt{\frac{x_i(x_{>i}+a_{>i})(a_i+1)_{x_i-1}(a_{>i}+1)_{x_{>i}}}{(x_i-1)!x_{>i}!}}
\end{align}
By using the relations
\begin{align*}
(a_i+1)_{x_i}&=\frac{(x_i+a_i)(a_i)_{x_i}}{a_i},\qquad \quad (a_i+1)_{x_i-1}=\frac{(a_i)_{x_i}}{a_i},\\
(a_{>i}+1)_{x_{>i}}&=\frac{(x_{>i}+a_{>i})(a_{>i})_{x_{>i}}}{a_{>i}},
\quad (a_{>i}+1)_{x_{>i}-1}=\frac{(a_{>i})_{x_{>i}}}{a_{>i}},\\
(-x_i+1)_{m-1-k}&=\frac{(-x_i)_{m-k}}{-x_i},\qquad \qquad \quad (-x_{>i}+1)_k=\frac{(-x_{>i})_{k+1}}{-x_{>i}},
\end{align*}
one arrives at
\begin{align*}
&\mathcal{B}(\bm{\lambda})^TX
 =\sqrt{\frac{(a_i)_{x_i}(a_{>i})_{x_{>i}}}{a_ia_{>i}x_i!x_{>i}!}}\n
&\qquad \qquad \qquad\times
\Bigl(-(x_i+a_i)(-x_i)_{m-1-k}(-x_{>i})_{k+1}+(x_{>i}+a_{>i})(-x_i)_{m-k}(-x_{>i})_k\Bigr),\\
&\quad -(x_i+a_i)(-x_i)_{m-1-k}(-x_{>i})_{k+1}\\
&\quad =-(x_i-m+k+1+a_i+m-k-1)(-x_i)_{m-1-k}(-x_{>i})_{k+1}\\
&\quad =(-x_i)_{m-k}(-x_{>i})_{k+1}-(a_i+m-k-1)(-x_i)_{m-1-k}(-x_{>i})_{k+1},  \tag{*}\\
&\quad (x_{>i}+a_{>i})(-x_i)_{m-k}(-x_{>i})_k=
 (x_{>i}-k+a_{>i}+k)(-x_i)_{m-k}(-x_{>i})_k\\
 & \quad =-(-x_i)_{m-k}(-x_{>i})_{k+1}+(a_{>i}+k)(-x_i)_{m-k}(-x_{>i})_k
\tag{**}.
\end{align*}
Cancellation of $(*)+(**)$ gives
\begin{align*}
&\mathcal{B}(\bm{\lambda})^TX\propto \phi_0(\bm{x};\bm{\lambda})\!
\times\!\Bigl((a_{>i}+k)(-x_i)_{m-k}(-x_{>i})_k\!-(a_i+m-k-1)(-x_i)_{m-1-k}(-x_{>i})_{k+1}\!\Bigr),
\end{align*}
and
\begin{align*}
&\sum_{k=0}^{m-1}(-1)^k\binom{m-1}{k}(a_{>i}+1+k)_{m-1-k}(a_i+m-k)_k\\
&\qquad \times\Bigl((a_{>i}+k)(-x_i)_{m-k}(-x_{>i})_k-(a_i+m-k-1)(-x_i)_{m-1-k}(-x_{>i})_{k+1} \Bigr)
\end{align*}
for the polynomial. For $k=0$ the first term in the parenthesis gives $(a_{>i})_m(-x_i)_m$, which is the 
$k=0$  term of the degree $m$ $\mathcal{P}_{m}^{(i)}(x_i,x_{>i};a_i,a_{>i})$.  
For $k-1$, the second term in the parenthesis gives 
\begin{equation*}
(-1)^k\binom{m-1}{k-1}(a_{>i}+k)_{m-k}(a_i+m-k)_k(-x_i)_{m-k}(-x_{>i})_k,\tag{*}
\end{equation*}
For $k$, the first term in the parenthesis turns out to be
\begin{equation*}
(-1)^k\binom{m-1}{k}(a_{>i}+k)_{m-k}(a_i+m-k)_k(-x_i)_{m-k}(-x_{>i})_k.\tag{**}
\end{equation*}
The sum of $(*)+(**)$ produces the $k$-th term of the degree $m$ $\mathcal{P}_{m}^{(i)}(x_i,x_{>i};a_i,a_{>i})$.
\hfill $\square$
\section*{Appendix B: Proof of the  eigenvalue formula, Theorem \ref{theo:eigens}}
\label{appendB}
\setcounter{equation}{0}
\renewcommand{\theequation}{B.\arabic{equation}}
Formula \eqref{eigj} is shown first. Let us evaluate 
$\bigl(\widetilde{\mathcal H}_j-\widetilde{\mathcal H}_i\bigr)\mathcal{P}_m^{(i)}$, $j<i$,
\begin{align*}
\bigl(\widetilde{\mathcal H}_j-\widetilde{\mathcal H}_i\bigr)\mathcal{P}_m^{(i)}
&=\bigl(\sum_{j\le k<i}x_k\bigr)(x_i+a_i)\bigl(\mathcal{P}_m^{(i)}-\mathcal{P}_m^{(i)}(x_i+1)\bigr)\\
&\quad \ +\bigl(\sum_{j\le k<i}x_k\bigr)(x_{>i}+a_{>i})\bigl(\mathcal{P}_m^{(i)}-\mathcal{P}_m^{(i)}(x_{>i}+1)\bigr)\\
&\quad \ +x_i\sum_{j\le k<i}(x_k+a_k)\bigl(\mathcal{P}_m^{(i)}-\mathcal{P}_m^{(i)}(x_{i}-1)\bigr)\\
&\quad \ +x_{>i}\sum_{j\le k<i}(x_k+a_k)\bigl(\mathcal{P}_m^{(i)}-\mathcal{P}_m^{(i)}(x_{>i}-1)\bigr),
\end{align*}
in which the default arguments of $\mathcal{P}_m^{(i)}$ are suppressed and those with changes are displayed.
\begin{align*}
\bigl(\widetilde{\mathcal H}_j-\widetilde{\mathcal H}_i\bigr)\mathcal{P}_m^{(i)}
&=\bigl(\sum_{j\le k<i}x_k\bigr)(x_i+x_{>i}+a_i+a_{>i})\mathcal{P}_m^{(i)} \tag{1}\\
&\ +\bigl(\sum_{j\le k<i}(x_k+a_k)\bigr)(x_i+x_{>i})\mathcal{P}_m^{(i)} \tag{2}\\
&\ -\bigl(\sum_{j\le k<i}x_k\bigr)\Bigl((x_i+a_i)\mathcal{P}_m^{(i)}(x_i+1)
+(x_{>i}+a_{>i})\mathcal{P}_m^{(i)}(x_{>i}+1)\Bigr) \tag{3}\\
&\ -\bigl(\sum_{j\le k<i}(x_k+a_k)\bigr)\Bigl(x_i\mathcal{P}_m^{(i)}(x_i-1)
+x_{>i}\mathcal{P}_m^{(i)}(x_{>i}-1)\Bigr) .  \tag{4}
\end{align*}
By using Theorem \ref{theo:rec}, (3) and (4) are simplified to
\begin{align*}
& -\bigl(\sum_{j\le k<i}x_k\bigr)(x_i+x_{>i}+a_i+a_{>i}+m) \mathcal{P}_{m}^{(i)} \tag{3'}\\
& -\bigl(\sum_{j\le k<i}(x_k+a_k)\bigr)(x_i+x_{>i}-m) \mathcal{P}_{m}^{(i)} \tag{4'}.
\end{align*}
Summing them and using \eqref{eigHi} gives 
\begin{equation*}
\widetilde{\mathcal H}_j\mathcal{P}_{m}^{(i)}=m(m+a_j+\cdots+a_n-1)\mathcal{P}_{m}^{(i)}.
\end{equation*}
The proof of the eigenvalue formula of $\widetilde{\mathcal H}_T$ \eqref{eigT} goes similarly.
As $\widetilde{\mathcal H}_T=\widetilde{\mathcal H}_0+\widetilde{\mathcal H}_1$ and the 
eigenvalue formula of $\widetilde{\mathcal H}_1$ is given above. By repeating the above argument for 
$\widetilde{\mathcal H}_0$ gives the result \eqref{eigT}.

\end{document}